\documentclass[12pt,a4paper]{article}
\usepackage{amsmath,amsthm,amssymb,latexsym}
\usepackage{enumerate}

\newtheorem{theorem}{Theorem}[section]
\newtheorem{lemma}[theorem]{Lemma}
\newtheorem{corollary}[theorem]{Corollary}
\newtheorem{proposition}[theorem]{Proposition}

\newtheorem{problem}[theorem]{Problem}
\newtheorem{remark}[theorem]{Remark}
\newtheorem{question}[theorem]{Question}

\renewcommand{\phi}{\varphi}

\begin{document}

\title{On classification of maximal nilpotent subsemigroups}
\author{Olexandr Ganyushkin and Volodymyr Mazorchuk}
\date{}
\maketitle

\begin{abstract}
We develop a general approach to the study of maximal
nilpotent subsemigroups of finite semigroups. This approach
can be used to recover many known classifications of maximal
nilpotent subsemigroups, in particular, for the symmetric
inverse semigroup, the symmetric semigroup, and the 
factor power of the symmetric group. We also apply this 
approach to obtain a classification of maximal nilpotent 
subsemigroups in the semigroup of binary relations and in 
certain $0$-simple finite semigroups.
\end{abstract}

\section{Introduction}\label{s1}

Let $S$ be a semigroup with the zero element $0$. A subsemigroup, 
$T\subset S$,
is called {\em nilpotent} provided that there exists $k\in\mathbb{N}$
such that $T^k=0$. The set of all nilpotent subsemigroups of $S$
is partially ordered with respect to inclusions, and the maximal
elements of this set are called {\em maximal nilpotent subsemigroups}
of $S$. The history of the problem to classify (or describe) all 
maximal nilpotent subsemigroups of a given semigroup goes back at
least to \cite{Gr}, where this problem was approached for 
$0$-simple semigroups using the graph theory. More recently the problem was 
studied for several classes of transformation semigroups in
\cite{BRR,GK1,GK2,GM1,GM2,Sh1,Sh2,St,Ts}. In all these cases a complete 
answer to the problem (that is a classification of all maximal 
nilpotent subsemigroups) was obtained. 
Although the technical details of the above papers are quite different,
one can single out the general method, used in all papers: 
the semigroups, which were studied, act (usually by partially defined
maps) on some set (not necessarily the one they are defined on) in 
some (natural) way; this action can be restricted to nilpotent 
subsemigroups, and, because of the nilpotency, defines a partial
order on the set; the maximal nilpotent subsemigroups can then be 
classified by the maximal partial orders (which are usually just the 
linear orders).

A natural question to ask now is: can this kind of approach be applied
to abstract semigroups? The aim of the present paper is to answer
this question positively. Under some technical conditions (which look
rather special but are satisfied, for example, by all finite inverse
semigroups) we develop a general machinery for classification of
maximal nilpotent subsemigroups of finite semigroups and show
that it can be applied to recover almost all the results listed 
above as well as to obtain many new results. 

The main idea of our approach is to find an appropriate set 
for the semigroup to ``act'' on. For this purpose we choose a 
special set of minimal idempotents and use the combinatorics of the 
egg-box diagram to define an ``action'' of the semigroup on
this set via certain matrices. The technical restrictions
on our semigroup we start with guarantee that the ``action''
is well-defined and faithfull. In this way we reduce the
study of the maximal nilpotent subsemigroups to the study of
partial (linear) orders on some sets of minimal idempotents.
The general classification then goes more or less along the 
same way as for transformation semigroups. 

The paper is organized as follows: In Section~\ref{s2}
we introduce all necessary notation and recall some
basic facts about nilpotent semigroups. In 
Section~\ref{s3} we define the radical of a semigroup
and study its properties. For a finite semigroup,
$S$, the radical $R(S)$ of $S$ is a nilpotent two-sided
ideal, which is contained in all maximal nilpotent
subsemigroups of $S$. This allows us to reduce 
the classification problem for the maximal nilpotent
subsemigroups of $S$ from $S$ to the Rees factor
$S/R(S)$. In Section~\ref{s4} we define and study the
properties of the so-called {\em minimal idempotents} of $S$. 
These are the ``points'' of the set the semigroup $S$
will ``act'' on. We suspect that for finite semigroups
the minimal idempotents are exactly the classical primitive ones. 
This is proved in Section~\ref{s4} for regular semigroups,
however, we do not know if this is true in the general
case. Section~\ref{s5} contains the main results
of our paper: in this section we give a classification
of all maximal nilpotent subsemigroups of finite
semigroups, satisfying some technical conditions
(conditions \eqref{cond1}-\eqref{cond3}). For the
semigroups with zero radical we even classify all
maximal subsemigroups among the nilpotent
subsemigroups of a fixed nilpotency class. Finally,
in Section~\ref{s6} we collected several applications
of our main results. Apart from the classical cases
which we recover (for example the cases of the
symmetric inverse semigroup and the semigroup of
all partial linear bijection on the finite-dimensional
vector space over a finite field), we also obtain
a classification of maximal nilpotent subsemigroups 
in the semigroup of all binary relations on a
finite set, and for several classes of finite $0$-simple
semigroups.
\vspace{0.2cm}

\noindent
{\bf Acknowledgments.}
The paper was written during the visit of the first author to 
Uppsala University, which was supported by The Swedish Institute.
The financial support of The Swedish Institute and the hospitality 
of Uppsala University are gratefully acknowledged. For the second 
author the research was partially supported by The Swedish 
Research Council.

\section{Notation and preliminaries}\label{s2}

Let $S$ be a semigroup. If the opposite is not explicitly stated, 
we assume that $S$ is finite an contains the zero element $0$. 

An element, $a\in S$, is called {\em nilpotent} provided that 
$a^n=0$ for some $n\in\mathbb{N}$. Analogously, $S$ is called
{\em nilpotent} provided that $S^n=0$ for some $n\in\mathbb{N}$.
The minimal $n\in \mathbb{N}$ such that $a^n=0$ (resp. $S^n=0$)
is called the {\em nilpotency class} of $a$ (resp. $S$). 

\begin{proposition}\label{prop1}{\rm (\cite[Fact~2.30,~page~179]{Ar})}
Let $S$ be a finite semigroup. Then the following statement
are equivalent:
\begin{enumerate}[(i)]
\item\label{prop1.1} $S$ is nilpotent.
\item\label{prop1.2} Every $a\in S$ is nilpotent.
\item\label{prop1.3} The only idempotent of $S$ is the zero element.
\end{enumerate}
\end{proposition}

We note that Proposition~\ref{prop1} fails for infinite semigroups
in general. For example, consider the semigroup
$T=\{(n,m):n,m\in\mathbb{N}, 0<n<m\}\cup\{0\}$ with the multiplication
\begin{displaymath}
(n,m)(k,l)=\begin{cases}
(n,l),&\text{if $m=k$},\\
0,&\text{otherwise}.
\end{cases}
\end{displaymath}
It is easy to see that every element in $T$ is nilpotent of class $2$, 
however, the semigroup $T$ is not nilpotent.

Recall the notation
\begin{displaymath}
S^1=\begin{cases}
S,& \text{if $S$ contains the identity element},\\
S\cup\{1\}, & \text{otherwise}.
\end{cases}
\end{displaymath}
The Green relations on $S$ will be denoted by $\mathcal{L}$,
$\mathcal{R}$, $\mathcal{H}$, $\mathcal{D}$, and $\mathcal{J}$.
For $a\in S$ we denote by $\mathcal{L}_a$ the $\mathcal{L}$-class
containing $a$ and analogously for other Green relations.
We will use the notation $(a)$ for the two-sided ideal $S^1aS^1$,
generated by $a\in S$. 

With the Green relations one associates partial quasi-orders
on $S$ in the following way (see e.g. \cite[II.1]{Gri}):
\begin{displaymath}
\begin{array}{rcl}
a\leq_{\mathcal{L}}b & \Leftrightarrow & S^1a\subset S^1b
\quad(\Leftrightarrow a=ub \text{ for some } u\in S^1),\\
a\leq_{\mathcal{R}}b & \Leftrightarrow  & aS^1\subset bS^1
\quad(\Leftrightarrow a=bu \text{ for some } u\in S^1),\\
a\leq_{\mathcal{H}}b & \Leftrightarrow  & a=bu=bv \text{ for some } 
u,v\in S^1,\\
a\leq_{\mathcal{J}}b & \Leftrightarrow  & (a)\subset (b)
\quad(\Leftrightarrow a=ubv \text{ for some } u,v\in S^1 ).
\end{array}
\end{displaymath}
Each of these quasi-orders induces in the natural way a partial order 
on the set of all equivalence classes with respect to the
corresponding Green relation. For finite semigroups the
relations $\mathcal{D}$ and $\mathcal{J}$ coincide and hence
$\leq_{\mathcal{J}}$ induces a partial order on the set of
all $\mathcal{D}$-classes. We will denote this order by $\leq$.
From the definition of $\leq_{\mathcal{J}}$ it follows that for every
$a\in S$ the principal ideal $(a)$ is the union of all those
$\mathcal{J}$-classes $J$ of $S$, which satisfy
$J\leq\mathcal{J}_a$ (for finite semigroups the same
is true for $\mathcal{D}$-classes).

We denote by $E(S)$ the set of all idempotents of $S$. The
restriction of the quasi-oder $\leq_{\mathcal{H}}$ to
$E(S)$ coincides with the so-called {\em natural} partial
order on $E(S)$, defined via: $f\preceq e$ if and only if
$ef=fe=f$, see \cite[page~195]{Gri}. The idempotents,
which are $0$-minimal with respect to this order,
are called {\em primitive}.

Assume now that $S$ is an arbitrary semigroup with the zero element
$0$ such that $\mathcal{D}=\mathcal{J}$.
A non-zero $\mathcal{D}$-class, $D$, of $S$ is called 
{\em minimal} provided that for any $\mathcal{D}$-class $D'$
the inequality $D'\leq D$ implies $D'=D$ or $D'=\{0\}$.

\begin{lemma}\label{lemma1}
Let $S$ be a finite semigroup, $D$ a minimal $\mathcal{D}$-class 
of $S$, and $a,b\in D$. Then either $ab=0$ or
$ab\in \mathcal{R}_a\cap \mathcal{L}_b$. In particular,
the set $D\cup \{0\}$ is a subsemigroup of $S$.
\end{lemma}

\begin{proof}
Assume $ab\neq 0$. Then the inclusion $(ab)\subset (a)$
and the minimality of $D=\mathcal{D}_a=\mathcal{D}_b$ imply 
that $ab\in D$. In particular, $b\mathcal{D}ab$, which 
yields $b=xaby$ for some $x,y\in S^1$. However, in this case
\begin{displaymath}
|S^1ab|\leq |S^1b|=|S^1xaby|\leq |S^1xab|\leq |S^1ab|,
\end{displaymath}
and hence the obvious inclusion $S^1ab\subset S^1b$ implies
$S^1ab=S^1b$, that is $ab\mathcal{L}b$. Analogously
one shows that $ab\mathcal{R}a$ and therefore 
$ab\in \mathcal{R}_a\cap \mathcal{L}_b$.
\end{proof}

\begin{remark}\label{rem1}
{\rm
The statement of Lemma~\ref{lemma1} is not true for infinite
semigroups in general. For example, in the bicyclic semigroup
\begin{equation}\label{eqbicyclic}
B=\langle a,b|ab=1\rangle=\{b^ka^m|k,m\geq 0\}
\end{equation}
we have
\begin{displaymath}
\mathcal{R}_{b^ka^m}=\{b^ka^n|n\geq 0\},\quad
\mathcal{L}_{b^ka^m}=\{b^na^m|n\geq 0\},
\end{displaymath}
and there is a unique $\mathcal{D}$-class. Hence in the semigroup
$B^0=B\cup\{0\}$ all non-zero elements form a unique
minimal $\mathcal{D}$-class. At the same time
$b^2=bb\not\in \mathcal{R}_b\cap\mathcal{L}_b=\{b\}$.
}
\end{remark}

\begin{corollary}\label{corollary2}
If a minimal $\mathcal{D}$-class, $D$, of a finite semigroup, $S$,
does not contain any idempotents (i.e. $D$ is not regular), then
the semigroup $D\cup\{0\}$ is nilpotent.
\end{corollary}

\begin{proof}
This is a direct corollary from Lemma~\ref{lemma1} and
Proposition~\ref{prop1}.
\end{proof}
 
\begin{lemma}\label{lemma2}
Let $S$ be an arbitrary semigroup such that 
$\mathcal{D}=\mathcal{J}$, and $D$ be a minimal 
$\mathcal{D}$-class. Then $(a)=D\cup \{0\}$ for
any $a\in D$.
\end{lemma}

\begin{proof}
Since $(x)=(y)$ for arbitrary $x,y\in D$, we have
$D\cup \{0\}\subset (a)$. On the other hand, for any
$0\neq b\in (a)$ we have $(b)\subset (a)$. However,
the minimality of $D$ implies the minimality of the 
principal ideal $(a)$ and hence $(a)=(b)$. This
means that $b\in \mathcal{D}_a=D$ and thus 
$(a)=D\cup \{0\}$.
\end{proof}

\section{The radical}\label{s3}

Let $S$ be an arbitrary semigroup. The set
\begin{displaymath}
R(S)=\{x:(x)\text{ is a nilpotent semigroup}\}
\end{displaymath}
will be called the {\em radical} of $S$.

\begin{lemma}\label{lemma3}
The radical $R(S)$ of an arbitrary semigroup, $S$, is
a (two-sided) ideal of $S$.
\end{lemma}

\begin{proof}
Let $x\in R(S)$ and $a,b\in S^1$. Then we have an obvious
inclusion $(axb)\subset (x)$. Since a subsemigroup of a
nilpotens semigroup is obviously nilpotent itself, we
get that $(axb)$ is nilpotent and hence $axb\in R(S)$.
\end{proof}

\begin{remark}\label{rem2}
{\em
If $S$ is inverse, then $R(S)=\{0\}$. Indeed, let
$a\in S$, $a\neq 0$. Then $(a)$ contains a non-zero
idempotent, namely $aa^{-1}\neq 0$, and hence $(a)$
is not nilpotent.
}
\end{remark}

Assume now that $S$ is such that $\mathcal{D}=\mathcal{J}$.
A $\mathcal{D}$-class, $D$, will be called {\em subminimal}
provided that for every $\mathcal{D}$-class $D'$ the
inequality $D'\leq D$ implies that either $D'=\{0\}$ or
$D'$ does not contain any idempotents (i.e. is not regular).
In particular, every minimal $\mathcal{D}$-class without
idempotents is subminimal.

\begin{lemma}\label{lemma4}
Let $S$ be finite. Then $R(S)$ coinsides with the union of
all subminimal $\mathcal{D}$-classes.
\end{lemma}

\begin{proof}
Since $R(S)$ is a two-sided ideal by Lemma~\ref{lemma3},
it is a union of $\mathcal{D}$-classes. Let $0\neq x\in R(S)$
and let $0\neq y\in S$ is such that $\mathcal{D}_y\leq \mathcal{D}_x$.
Then $(y)\subset (x)$, in particular $(y)$ is nilpotent, implying
that the only idempotent of $(y)$ is $0$. Hence $y$ is not
an idempotent and $\mathcal{D}_x$ is subminimal.

Let now $D$ be a subminimal $\mathcal{D}$-class and $x\in D$.
The ideal $(x)$ is the union of all $\mathcal{D}$-classes $D'$ such
that $D'\leq D$ and hence $(x)$ does not contain non-zero idempotents
by the subminimality of $D$. Hence $(x)$ is nilpotent and $x\in R(S)$.
Since $x\in D$ was arbitrary, we even obtain $D\subset R(S)$.
\end{proof}

\begin{corollary}\label{cornew}
The radical of a finite semigroup is a nilpotent semigroup.
\end{corollary}

\begin{proof}
$R(S)$ is a semigroup by Lemma~\ref{lemma3}.
From Lemma~\ref{lemma4} it follows that the only idempotent
of $R(S)$ is zero. Hence $R(S)$ is nilpotent by Proposition~\ref{prop1}.
\end{proof}

\begin{remark}\label{remnilp}
{\rm 
Corollary~\ref{cornew} is not true for infinite semigroups in the
general case. For $n\in\mathbb{N}$ let $T_n$ denote the Rees
factor of $(\mathbb{N},+)$ modulo the ideal $\{n,n+1,\dots\}$.
Let $S$  be the disjoint union of all $T_n$, $n\in\mathbb{N}$,
with the common zero element. Then every element of $S$ is nilpotent
by definition and it is easy to see that $S=R(S)$. However, the 
nilpotency classes of the elements in $S$ are not bounded and
hence $S$ itself is not nilpotent.
}
\end{remark}

\begin{lemma}\label{lemma5}
Let $S$ be arbitrary. Assume that
$R(S)$ is a nilpotent semigroup of class $k$. Let
$T\subset S$ be a nilpotent subsemigroup of $S$ of class $m$.
Then the subsemigroup $\langle R(S),T\rangle$ coinsides with
$R(S)\cup T$ and is a nilpotent subsemigroup of $S$ of
class at most $mk$.
\end{lemma}

\begin{proof}
The equality $\langle R(S),T\rangle=R(S)\cup T$ follows from the
fact that $R(S)$ is an ideal (Lemma~\ref{lemma3}). Let
$a_i\in R(S)\cup T$, $i=1,\dots,mk$. If at least 
$k$ elements out of $a_1,\dots,a_{mk}$ belong to $R(S)$,
the product $a_1\cdots a_{mk}$ reduces to a product of 
$k$ elements from $R(S)$ and hence equals $0$ since 
$R(S)$ is nilpotent of class $k$. On the other
hand, if the number of elements from $R(S)$ among
$a_1,\dots,a_{mk}$ is smaller than $k$, the equality
$k-1+k(m-1)=km-1$ implies that the product 
$a_1\cdots a_{mk}$ contains a subproduct, consisting
of at least $m$ consecutive factors from $T$. Hence  
$a_1\cdots a_{mk}=0$ since $T$ is nilpotent of class $m$.
\end{proof}

\begin{corollary}\label{cor3}
Let $S$ be an arbitrary semigroup such that $R(S)$ is
nilpotent. Then every maximal nilpotent subsemigroup of $S$ contains
$R(S)$. In particular, the canonical epimorphism $\varphi$
from $S$ to the Rees factor $S/R(S)$ induces a 
bijection between the maximal nilpotent subsemigroups of $S$
and the maximal nilpotent subsemigroups of $S/R(S)$.
\end{corollary}

\begin{proof}
Follows immediately from Lemma~\ref{lemma5} and the fact that
$\varphi$ sends nilpotent semigroups to nilpotent semigroups. 
\end{proof}

\begin{remark}\label{rem3}
{\em 
By Corollary~\ref{cor3}, the radical $R(S)$ of a semigroup, $S$, 
is contained in the intersection of all maximal nilpotent 
subsemigroups of $S$ provided that $R(S)$ is nilpotent. 
We will later see that in many cases $R(S)$ coincides
with the  intersection of all maximal nilpotent subsemigroups of
$S$. However, this is not true in the general case, as follows
from the following example:

Consider the following subsemigroup of the semigroup of all
partial transformations of the set $\{1,2,3,4\}$:
\begin{displaymath}
S=\left\{
0,
\left(\begin{array}{cc}1&3\\1&1\end{array}\right),
\left(\begin{array}{cc}1&3\\2&2\end{array}\right),
\left(\begin{array}{cc}1&3\\3&3\end{array}\right),
\left(\begin{array}{cc}1&3\\4&4\end{array}\right),
\overline{1},\overline{2},\overline{3},\overline{4}
\right\},
\end{displaymath}
where $\overline{k}$ denotes the constant transformation
$\left(\begin{array}{cccc}1&2&3&3\\k&k&k&k\end{array}\right)$.
It is easy to see that $S$ is indeed a semigroup. Observe
that $\overline{1}$ is an inverse element to both
$\left(\begin{array}{cc}1&3\\2&2\end{array}\right)$ and
$\left(\begin{array}{cc}1&3\\4&4\end{array}\right)$, and that all
other elements of $S$ are idempotents. Hence $S$ is regular,
which means that for every $0\neq x\in S$ the ideal
$(x)$ contains a non-zero idempotent (for example
$\overline{1}$ if 
$x=\left(\begin{array}{cc}1&3\\2&2\end{array}\right)$ or
$x=\left(\begin{array}{cc}1&3\\4&4\end{array}\right)$, and
$x$ itself in the other cases). This implies
$R(S)=0$. On the other hand, $S$ contains the unique maximal
nilpotent subsemigroup, namely
\begin{displaymath}
T=\left\{
0,
\left(\begin{array}{cc}1&3\\2&2\end{array}\right),
\left(\begin{array}{cc}1&3\\4&4\end{array}\right)
\right\}.
\end{displaymath}
}
\end{remark}

\section{Minimal idempotents}\label{s4}

From Corollary~\ref{cor3} it follows that for the classification
of  maximal nilpotent subsemigroups of a finite semigroup,
$S$, it is enough to solve this problem for the Rees quotient
$S/R(S)$. From Lemma~\ref{lemma4} it follows that in this
case every minimal $\mathcal{D}$-class of $S/R(S)$ contains idempotents.
The idempotents of $S$, which are contained in the minimal $\mathcal{D}$-classes of the semigroup $S/R(S)$ will
be called {\em minimal}. 

From now on we assume that 
$S$ is finite. Since $R(S/R(S))=\{0\}$, until the end of this
section we simply assume that  $R(S)=\{0\}$.

\begin{lemma}\label{lemma6}
Let $e,f\in E(S)$ be minimal such that $\mathcal{D}_e\neq\mathcal{D}_f$. 
Then $exf=0$ for any $x\in S^1$.
\end{lemma}

\begin{proof}
$exf\in (e)\cap (f)$. Since $\mathcal{D}_e\neq\mathcal{D}_f$, we get
$(e)\neq (f)$ and hence $(e)\cap (f)=\{0\}$ by the
minimality of $e$ and $f$.
\end{proof}

\begin{corollary}\label{cor4}
Let $e,f\in E(S)$ be minimal and $x\in S^1$ be such that
$exf\neq 0$. Then $\mathcal{D}_e=\mathcal{D}_f$ and
$exf\in \mathcal{D}_e$.
\end{corollary}

\begin{proof}
The equality $\mathcal{D}_e=\mathcal{D}_f$ follows from
Lemma~\ref{lemma6}. The inclusion
$(exf)=S^1exfS^1\subset S^1eS^1=(e)$ and the minimality of 
the ideal $(e)$ implies $(exf)=(e)$, and hence
$e\mathcal{D}exf$. This completes the proof.
\end{proof}

Two idempotents, $e,f\in E(S)$ will be called {\em orthogonal}
provided that $ef=fe=0$.

\begin{proposition}\label{prop2}
Let $e,f\in E(S)$ be minimal and $e\neq f$. 
Then $e$ and $f$ are orthogonal
if and only if they commute.
\end{proposition}

\begin{proof}
The necessity is obvious. To prove the sufficiency let us assume
that $ef=fe=a$. Then $a^2=efef=eeff=ef=a$ and hence $a\in E(S)$.
Assume that $a\neq 0$. Then $a=ef=eeff=eaf\neq 0$ and 
Corollary~\ref{cor4} implies $e\mathcal{D}f$ and 
$e\mathcal{D}a$. The minimality of $\mathcal{D}_e$ even gives
$(e)=(f)=(a)$. Moreover, $ea=ae=a$ and $af=fa=a$. 
Hence 
\begin{equation}\label{eqprop2.1}
S^1a=S^1ae\subset S^1e, \quad S^1a\subset S^1f,\quad
aS^1\subset S^1e, \quad aS^1\subset fS^1.
\end{equation}
On the other hand, $(e)=(a)$ implies that $e=xay$ for some
$x,y\in S^1$. Hence $S^1e=S^1xay$, implying
\begin{displaymath}
|S^1a|\geq |S^1xa|\geq |S^1xay|=|S^1e|.
\end{displaymath}
Together with \eqref{eqprop2.1} this gives $S^1a=S^1e$, that
is $a\mathcal{L}e$. Analogously we get 
$a\mathcal{L}f$, $a\mathcal{R}e$, and $a\mathcal{R}f$. This
implies $e\mathcal{H}f$ and thus $e=f$, a contradiction.
Therefore $a=0$ and the proof is complete.
\end{proof}

\begin{lemma}\label{lemma7}
Let $e\in E(S)$ be minimal and $f\in E(S)$, $f\neq 0$. Then
$f\preceq e$ implies $f\in\mathcal{D}_e$, in particular,
$f$ is minimal.
\end{lemma}

\begin{proof}
$f\preceq e$ implies $fe=ef=f$. Hence $efe=f\neq 0$ and
$f\in \mathcal{D}_e$ by Corollary~\ref{cor4}.
\end{proof}

We would like to finish this section with the study of
the relation between the primitive and the minimal
idempotents. 

\begin{remark}\label{rem03}
{\rm
For infinite semigroups, 
even inverse ones, minimal idempotents do not have to be 
primitive. Indeed, the bicyclic semigroup $B$, considered
in Remark~\ref{rem1} (see 
\eqref{eqbicyclic}), is a simple inverse semigroup.
Hence all non-zero elements of the semigroup
$B^0=B\cup\{0\}$ form a unique $\mathcal{D}$-class, which
is obviously minimal. Therefore all idempotents of 
$B$ are minimal in $B^0$. On the other hand, the idempotents
of $B$ form an infinite decreasing chain with respect to $\preceq$
and hence $B^0$ does not contain any primitive idempotent.
}
\end{remark}

For finite semigroups the situation is completely different:

\begin{theorem}\label{theorem1}
Let $S$ be a finite semigroup. Then every minimal idempotent 
of $S$ is primitive.
\end{theorem}

\begin{proof}
Without loss of generality we may assume $R(S)=\{0\}$.
Let $e,f\in E(S)$ be such that $e$ is minimal and
$ef=fe=f\neq 0$. Then $f\in \mathcal{D}_e$ by Lemma~\ref{lemma7}.
Then the minimality of $\mathcal{D}_e$ and Lemma~\ref{lemma1}
imply
\begin{displaymath}
f\in \mathcal{R}_e\cap \mathcal{L}_f\quad
\text{and}\quad
f\in \mathcal{R}_f\cap \mathcal{L}_e.
\end{displaymath}
Hence $f\in \mathcal{R}_e\cap \mathcal{L}_e=\mathcal{H}_e$ and
therefore $f=e$. This implies that $e$ is primitive.
\end{proof}

\begin{theorem}\label{theorem2}
Let $S$ be a regular semigroup such that $\mathcal{D}=\mathcal{J}$. 
Then every primitive idempotent of $S$ is minimal.
\end{theorem}

\begin{proof}
Assume that $e\in E(S)$ is primitive but not minimal. For the right
ideal $eS^1$ we have two possible cases.

{\bf Case 1.} Assume that $eS^1$ is not a minimal right ideal.
Then $eS^1$ properly contains some non-zero ideal $fS^1$. Since
$S$ is regular, we can even assume that $f$ is an idempotent.
Hence there exists $a\in S^1$ such that $ea=f=f^2=eaea\neq 0$,
moreover $eaS^1\subset eS^1$. From the equalities
$eae\cdot eae=eaeae=eae$ it follows that $eae$ is an 
idempotent too. Moreover, $eae\neq 0$ since 
$eae\cdot a=(ea)^2\neq 0$. Obviously $e\cdot eae=eae\cdot e=eae$,
that is $eae\preceq e$. Moreover, 
\begin{displaymath}
eaeS^1\subset eaS^1=fS^1\neq eS^1,
\end{displaymath}
and hence $eae\neq e$. This implies $0\prec eae\prec e$, which
contradicts our assumption that $e$ is primitive.

{\bf Case 2.} Assume that $eS^1$ is a minimal right ideal.
Then $eS^1=\mathcal{R}_e\cup\{0\}$ and thus
\begin{displaymath}
(e)=S^1eS^1=\cup_{x\in\mathcal{R}_e}S^1x.
\end{displaymath}
Since $e$ is not minimal, $(e)$ contains a $\mathcal{D}$-class,
say $D$, such that $D< \mathcal{D}_e$. Since $R$ is regular,
$D$ contains an idempotent, say $f$, that is $f=yx$ for some
$y\in S^1$ and $x\in \mathcal{R}_e$. From $x\mathcal{R}e$ we
have $x=ea$ and $e=xb$ for some $a,b\in S^1$. Hence
$f=yx=yea=yxba$. Since $(yx)\supset (yxb)\supset(yxba)=(yx)$,
we obtain
\begin{equation}\label{eqxxx1}
(f)=(yx)=(yxb)=(ye). 
\end{equation}
Obviously, $S^1ye\subset S^1e$. However, \eqref{eqxxx1}
implies that $ye\in \mathcal{D}_f=D\neq \mathcal{D}_e$
and hence $S^1ye\neq S^1e$. Thus $0\neq S^1ye\subsetneqq S^1e$
and the left ideal $S^1e$ is not minimal. Analogously to 
Case 1 one now shows that $e$ is not primitive, a contradiction.
This completes the proof.
\end{proof}

\begin{corollary}\label{cor5}
Let $S$ be finite and regular. Then $e\in E(S)$ is minimal if and
only if it is primitive.
\end{corollary}

\begin{question}\label{quest1}
Is the statement of Theorem~\ref{theorem2} true for
all finite semigroups? 
\end{question}

\section{Maximal nilpotent subsemigroups}\label{s5}

As in the previous section, we assume that $S$ is finite
and contains $0$. 

\begin{lemma}\label{lemma7.1}
Let $x\not\in R(S)$. Then there exist $e,f\in E(S)$ such that
$exf\neq 0$.
\end{lemma}

\begin{proof}
Since the ideal $(x)$ is not nilpotent, it must contain a non-zero
idempotent, say $g=axb$, where $a,b\in S^1$. In particular $g^k=g$
for all $k\in\mathbb{N}$.

As $S$ is finite, there exist $n,m\in\mathbb{N}$ such that
$e=(xba)^m$ and $f=(bax)^n$ are idempotents. However, we have
\begin{displaymath}
g=g^{n+m+1}=(axb)^{n+m+1}=a\cdot (xba)^m\cdot x\cdot (bax)^n\cdot b=
aexfb\neq 0.
\end{displaymath}
Therefore $exf\neq 0$.
\end{proof}

\begin{lemma}\label{lemma8}
Assume that $R(S)=0$ and  $0\neq x\in S$. Then there exist 
$e,f\in E(S)$ such that $exf\neq 0$, moreover, $e$ and $f$ can
be chosen in the same minimal $\mathcal{D}$-class. 
\end{lemma}

\begin{proof}
If $R(S)=0$, then the idempotent $g$ in the proof of Lemma~\ref{lemma7.1}
can be chosen from a minimal $\mathcal{D}$-class. Moreover, the 
numbers $m$ and $n$ from the definition of $e$ and $f$ can be
chosen $\geq 2$. In this case
\begin{displaymath}
e=xba\cdot xba\cdot (xba)^{m-2}=xb\cdot g\cdot a(xba)^{m-2}
\end{displaymath}
Hence $e\in (g)$ and $\mathcal{D}_e\leq \mathcal{D}_g$. The minimality
of $\mathcal{D}_g$ thus yields $\mathcal{D}_e= \mathcal{D}_g$.

Analogously one shows that $\mathcal{D}_f= \mathcal{D}_g$ and
the statement is proved.
\end{proof}

From now on we assume that $R(S)=0$. Assume that we are given 
a non-empty set, $\mathcal{M}$, of minimal idempotents of
$S$, which satisfies the following conditions:
\begin{enumerate}[(I)]
\item\label{cond1} the elements of $\mathcal{M}$ commute;
\item\label{cond2} for every $0\neq x\in S$ there exist
$e,f\in \mathcal{M}$ such that $exf\neq 0$;
\item\label{cond3} for arbitrary $x,y\in S$ and arbitrary
$e,f\in \mathcal{M}$ there is the following inclusion:
\begin{displaymath}
\{exgyf:g\in \mathcal{M}\}\subset \{0,exyf\},
\end{displaymath}
moreover, $\{exgyf:g\in \mathcal{M}\}\neq 0$ if 
$exyf\neq 0$.
\end{enumerate}
Consider the set $\mathrm{Bin}(S,\mathcal{M})$, which consists
of all $\mathcal{M}\times \mathcal{M}$-matrices 
$A=(a_{e,f})_{e,f\in \mathcal{M}}$,
where $a_{e,f}\in eSf$. Let
$A\in \mathrm{Bin}(S,\mathcal{M})$. If there exists
some $x\in S$ such that $a_{e,f}=exf$ for all $f\in \mathcal{M}$
(resp. for all $e\in \mathcal{M}$), then the corresponding
row (resp. column) of $A$ will be called {\em stable} and denoted
by $e_x$ (resp. ${}_xf$). Let now $e_x$ be a stable row of $A$
and ${}_yf$ be a stable column of some $B\in 
\mathrm{Bin}(S,\mathcal{M})$. If there exists $g\in\mathcal{M}$
such that $exgyf\neq 0$, we will say that the {\em product}
$e_x\cdot {}_yf$ of $e_x$ and ${}_yf$ equals $exgyf$. Otherwise we
set $e_x\cdot {}_yf=0$. From the equality $exgyf=exg\cdot gyf$ it 
follows that the value of $e_x\cdot {}_yf$ does not depend on the 
choice of $x$ and $y$, and from the condition \eqref{cond3} it follows
that $e_x\cdot {}_yf$ does not depend on the choice of $g$. Hence
the product of a stable row and a stable column is well-defined.
If either the $e$-th row of $A$ is not stable or the $f$-th column of
$B$ is not stable  (or both), we set that their product is equal
to  $0$. In this way we define the product of any two matrices
$A,B\in \mathrm{Bin}(S,\mathcal{M})$.

Consider the map
\begin{displaymath}
\begin{array}{cccc}
\psi: & S & \rightarrow & \mathrm{Bin}(S,\mathcal{M})\\
& x & \mapsto & A_x=(exf)_{e,f\in\mathcal{M}}.
\end{array}
\end{displaymath}

\begin{proposition}\label{prop4}
The map $\psi$ is a homomorphism.
\end{proposition}

\begin{proof}
First we remark that all rows and columns of $A_x$ are stable
by definition. Consider now some element $exyf$ of the matrix
$A_{xy}$. If $exyf=0$, from the condition \eqref{cond3} we
obtain that $exgyf=0$ for any $g\in \mathcal{M}$. Hence
the product of the row $e_x$ of $A_{x}$ with the column
${}_yf$ of $A_y$ equals $0$ as well. If $exyf\neq 0$, from 
the condition \eqref{cond3} we obtain that $exgyf=exyf$ 
for some $g\in \mathcal{M}$. Hence from the definition of
the product in $\mathrm{Bin}(S,\mathcal{M})$ it follows that
in this case $e_x\cdot {}_yf=exyf$. This shows that 
$A_{xy}=A_x\cdot A_y$ and completes the proof.
\end{proof}

Conditions \eqref{cond1}-\eqref{cond3} might look exotic, however,
they can be satisfied for many semigroups and even classes of 
semigroups. For instance, we have the following statement:

\begin{theorem}\label{theorem3}
The set of all minimal idempotents of a finite inverse semigroup
satisfies the conditions \eqref{cond1}-\eqref{cond3}.
\end{theorem}

\begin{proof}
Let $S$ be a finite inverse semigroup, and $\mathcal{M}$ be the
set of all minimal idempotents of $S$. Then \eqref{cond1} is
satisfied as all idempotents of an inverse semigroup commute.
Since $R(S)=0$ by Remark~\ref{rem2}, the condition 
\eqref{cond2} follows from Lemma~\ref{lemma7.1} and Lemma~\ref{lemma8}.
So, we are left to prove that the condition \eqref{cond3} is
also satisfied. By the Wagner-Preston Theorem, we can consider $S$
as a subsemigroup of the inverse semigroup $\mathcal{IS}(S)$ of
all partial bijections on $S$. We will use the standard notation
$\mathrm{dom}$ and $\mathrm{im}$ for the {\em domain} and the
{\em range} of partial maps. By Theorem~\ref{theorem2}, the 
minimal and the primitive idempotents of $S$ coincide. 

Let $x,y\in S$ and $e,f$ be minimal idempotents. Then for any 
minimal idempotent $g$ we have either $\mathrm{dom}(g)\cap
\mathrm{im}(ex)=\varnothing$ or $\mathrm{dom}(g)=\mathrm{im}(ex)$.
Indeed, let $\mathrm{dom}(g)\cap\mathrm{im}(ex)\neq\varnothing$.
Then $exg\neq 0$ and the non-zero idempotent
$h=(exg)\cdot(exg)^{-1}$ satisfies $f\preceq e$. Since $e$
is primitive, we have $h=e$, which is possible only if
$\mathrm{dom}(g)\supset\mathrm{im}(ex)$. Analogously one proves
$\mathrm{dom}(e)\supset\mathrm{im}(gx^{-1})$, which implies
$\mathrm{dom}(g)\subset\mathrm{im}(ex)$. Hence
$\mathrm{dom}(g)=\mathrm{im}(ex)$.

If $\mathrm{dom}(g)\cap\mathrm{im}(ex)=\varnothing$, we have
$exgyf=0$. If $\mathrm{dom}(g)=\mathrm{im}(ex)$, we have
$exg=ex$ and hence $exgyf=exyf$. Therefore
\begin{displaymath}
\{exgyf:g\in\mathcal{M}\}\subset \{0,exyf\}.
\end{displaymath}
On the other hand, if $exyf\neq 0$, then
$ex\neq 0$ and $g=(ex)^{-1}\cdot ex$ is a non-zero idempotent.
From the inequality $h<g$ it follows that 
$exh\cdot (exh)^{-1}<e$ and hence $h$ is primitive. But
$exhyf=ex\cdot (ex)^{-1}\cdot ex\cdot yf=exyf$. Hence
$exyf\neq 0$ implies $\{exgyf:g\in\mathcal{M}\}\neq \{0\}$,
which gives us \eqref{cond3} and completes the proof.
\end{proof}

Let us now assume that we are given a finite semigroup, $S$,
such that $R(S)=0$ and such that there exists a non-empty set,
$\mathcal{M}$, of minimal idempotents of $S$, satisfying the
conditions \eqref{cond1}-\eqref{cond3}. Denote by
$\mathrm{Nil}(S)$ the set of all nilpotent subsemigroups of
$S$, and by $\mathrm{Ord}(\mathcal{M})$ the set of all 
(strict) partial orders on $\mathcal{M}$. Both
$\mathrm{Nil}(S)$ and $\mathrm{Ord}(\mathcal{M})$ are partially
ordered with respect to inclusions in the natural way.

For every $\rho\in \mathrm{Ord}(\mathcal{M})$ set
\begin{equation}\label{eqord}
\mathrm{Mon}(\rho)=\{x\in S: exf\neq 0\text{ implies }
(e,f)\in\rho\text{ for all }e,f\in \mathcal{M}\}.
\end{equation}
For every $T\in \mathrm{Nil}(S)$ set
\begin{equation}\label{eqnil}
\rho_T=\{(e,f)|exf\neq 0\text{ for some }x\in T\}
\subset \mathcal{M}\times \mathcal{M}.
\end{equation}

\begin{proposition}\label{prop5}
\begin{enumerate}[(a)]
\item\label{prop5.1} The map $\rho\mapsto \mathrm{Mon}(\rho)$
is a homomorphism from the poset
$\mathrm{Ord}(\mathcal{M})$ to the poset $\mathrm{Nil}(S)$.
\item\label{prop5.2} The map $T\mapsto \rho_T$
is a homomorphism from the poset
$\mathrm{Nil}(S)$ to the poset $\mathrm{Ord}(\mathcal{M})$.
\end{enumerate}
\end{proposition}

\begin{proof}
Let us prove \eqref{prop5.1}. To show that $\mathrm{Mon}(\rho)$
is a subsemigroup it is enough to show that it is closed with respect
to multiplication. Let $x,y\in \mathrm{Mon}(\rho)$ and
$exyf\neq 0$ for some $e,f\in \mathcal{M}$. Then from \eqref{cond3}
we have that there exists $g\in \mathcal{M}$ such that
$exgyf=exyf$. This implies  $exg\neq 0$ and $gyf\neq 0$ and hence
$e\rho g$ and $g\rho f$. The transitivity of $\rho$ yields
$e\rho f$ and thus $xy\in \mathrm{Mon}(\rho)$.

Let $n=|\mathcal{M}|$ and $x_1,\dots,x_n\in \mathrm{Mon}(\rho)$
be arbitrary. Assume that the exist $e,f\in \mathcal{M}$ such that
$ex_1\cdots x_nf\neq 0$. From \eqref{cond3} we obtain the existence
of $g_1,\dots,g_{n-1}\in\mathcal{M}$ such that
$ex_1g_1x_2\cdots x_{n-1}g_{n-1}x_nf\neq 0$. Hence
$ex_1g_1\neq 0$, $g_1x_2g_2\neq 0$,\dots,
$g_{n-1}x_nf\neq 0$, implying $e\rho g_1$, $g_1\rho g_2$,\dots,
$g_{n-1}\rho f$. This gives us a chain of cardinality $n+1$ in
$\mathcal{M}$, which is not possible, a contradiction. Hence
$ex_1\cdots x_nf= 0$ for all $e,f\in \mathcal{M}$ and therefore
$x_1\cdots x_n=0$ by \eqref{cond2}. This implies that
$\mathrm{Mon}(\rho)$ is a nilpotent subsemigroup of $S$ of
nilpotency class at most $n$. 

The implication $\rho_1\subset\rho_2$ implies
$\mathrm{Mon}(\rho_1)\subset \mathrm{Mon}(\rho_2)$ is obvious.
This proves \eqref{prop5.1}.

Let us now prove \eqref{prop5.2}. Let $(e,f)$ and $(f,g)$
belong to $\rho_T$. Then there exist $x,y\in S$ such that
$exf\neq 0$ and $fyg\neq 0$. Lemma~\ref{lemma6} implies that 
$\mathcal{D}_e=\mathcal{D}_f=\mathcal{D}_g$. The
equalities $e\cdot exf=exf\cdot f=exf$ and Lemma~\ref{lemma1}
imply that $exf\in \mathcal{R}_e\cap\mathcal{L}_f$.
Analogously $fyg\in \mathcal{R}_f\cap\mathcal{L}_g$.
Since $f$ is an idempotent, $exf\in \mathcal{L}_f$ and
$fyg\in \mathcal{R}_f$,  from \cite[Proposition~2.4]{Gri}
we obtain that $exf\cdot fyg=exfyg\in \mathcal{R}_{exf}\cap
\mathcal{L}_{fyg}=\mathcal{R}_e\cap\mathcal{L}_{g}$.
In particular, $exfyg\neq 0$. From \eqref{cond3} we get that
$exyg\neq 0$. As $xy\in T$, we obtain $(e,g)\in\rho_T$ and hence
$\rho_T$ is transitive.

Assume that for some $e,f\in \mathcal{M}$ we have both
$(e,f)\in\rho_T$ and $(f,e)\in\rho_T$. Then, analogously to the
arguments above, we find $x,y\in T$ such that 
$exf\neq 0$, $fye\neq 0$ and $exye\neq 0$. We have
$exye\in \mathcal{H}_e$ by Lemma~\ref{lemma1}. Moreover, 
$\mathcal{H}_e$ contains an idempotent, which means that
$\mathcal{H}_e$ is a group. Hence $(exye)^m=e(xy)^me\in 
\mathcal{H}_e$ for all $m\in\mathbb{N}$. This implies
$(xy)^m\neq 0$ for all $m\in\mathbb{N}$, which contradicts the
nilpotency of $T$. This proves that $\rho_T$ is anti-symmetric
ands hence $\rho_t\in \mathrm{Ord}(\mathcal{M})$.

That the implication $T_1\subset T_2$ implies
$\rho_{T_1}\subset \rho_{T_2}$ is obvious and the proof is 
complete.
\end{proof}

\begin{proposition}\label{prop6}
Let $T\in\mathrm{Nil}(S)$ and $\rho\in \mathrm{Ord}(\mathcal{M})$
be arbitrary. Then
\begin{enumerate}[(a)]
\item\label{prop6.1} $T\subset \mathrm{Mon}(\rho_T)$,
$\rho_{\mathrm{Mon}(\rho)}\subset \rho$;
\item\label{prop6.2} $\mathrm{Mon}(\rho_{\mathrm{Mon}(\rho)})=
\mathrm{Mon}(\rho)$;
\item\label{prop6.3} $\rho_{\mathrm{Mon}(\rho_T)}=\rho_T$.
\end{enumerate}
\end{proposition}

\begin{proof}
The statement \eqref{prop6.1} is obvious.

As $\rho_{\mathrm{Mon}(\rho)}\subset \rho$, from
Proposition~\ref{prop5}\eqref{prop5.1} it follows that
$\mathrm{Mon}(\rho_{\mathrm{Mon}(\rho)})\subset
\mathrm{Mon}(\rho)$. Let now $x\in \mathrm{Mon}(\rho)$.
From the definition \eqref{eqnil} we have that for all
$e,f\in\mathcal{M}$ from $exf\neq 0$ it follows that
$(e,f)\in \rho_{\mathrm{Mon}(\rho)}$. According to
\eqref{eqord}, this implies $x\in 
\mathrm{Mon}(\rho_{\mathrm{Mon}(\rho)})$, which proves
$\mathrm{Mon}(\rho)\subset 
\mathrm{Mon}(\rho_{\mathrm{Mon}(\rho)})$. The statement
\eqref{prop6.2} follows.

As $T\subset \mathrm{Mon}(\rho_T)$, from
Proposition~\ref{prop5}\eqref{prop5.2} it follows that
$\rho_T\subset \rho_{\mathrm{Mon}(\rho_T)}$. Let now 
$(e,f)\in \rho_{\mathrm{Mon}(\rho_T)}$.
From the definition \eqref{eqnil} it follows that there
exists $x\in \mathrm{Mon}(\rho_T)$ such that $exf\neq 0$.
This and the definition of $\mathrm{Mon}(\rho_T)$ implies 
$(e,f)\in\rho_T$. Hence $\rho_{\mathrm{Mon}(\rho_T)}\subset
\rho_T$ and the statement \eqref{prop6.3} follows. This
completes the proof.
\end{proof}

Recall that, according to \cite{Co}, a pair of maps, $\alpha:P\to Q$
and $\beta:Q\to P$, defines a {\em Galois correspondence} between the posets
$P$ and $Q$ if it satisfies the following conditions:
\begin{enumerate}
\item $\alpha$ and $\beta$ are antihomomorphisms of the partially
ordered sets, that is $p_1\leq p_2$ implies $\alpha(p_1)\geq
\alpha(p_2)$ and $q_1\leq q_2$ implies $\beta(q_1)\geq \beta(q_2)$;
\item $\alpha\beta(p)\geq p$ and $\beta\alpha(q)\geq q$ for all
$p\in P$ and $q\in Q$;
\item $\beta\alpha\beta(p)=\beta(p)$ and
$\alpha\beta\alpha(q)=\alpha(q)$ for all
$p\in P$ and $q\in Q$.
\end{enumerate}

Denote by $\mathrm{Ord}(\mathcal{M})^*$ the set 
$\mathrm{Ord}(\mathcal{M})$ with the order, which is
the opposite to the inclusion order.

\begin{theorem}\label{theorem3.1}
The pair of maps $T\mapsto\rho_T$ and
$\rho\mapsto\mathrm{Mon}(\rho)$ defines a Galois
correspondence between the posets
$\mathrm{Nil}(S)$ and $\mathrm{Ord}(\mathcal{M})^*$.
\end{theorem}

\begin{proof}
From Proposition~\ref{prop5} it follows that 
the maps $T\mapsto\rho_T$ and
$\rho\mapsto\mathrm{Mon}(\rho)$ are antihomomorphisms
between the partially ordered sets
$\mathrm{Nil}(S)$ and $\mathrm{Ord}(\mathcal{M})^*$.
The rest follows from Proposition~\ref{prop6}.
\end{proof}

Since the restriction of $\mathcal{D}$ to an arbitrary
subset, $P\subset S$, defines on $P$ some equivalence
relation, we can consider the decomposition
$P=P_1\cup\dots\cup P_k$ of $P$ into  a disjoint union
of equivalence classes with respect to this relation.
The main result of the present paper is the following theorem:

\begin{theorem}\label{theorem4}
Let $S$ be a finite semigroup with radical $R(S)$, and
$\mathcal{M}$ be a non-empty set of minimal idempotents of
$S/R(S)$, which satisfies \eqref{cond1}-\eqref{cond3}. Let
$\mathcal{M}=\mathcal{M}_1\cup\dots\cup \mathcal{M}_k$ be
a decomposition of $\mathcal{M}$ into a disjoint union of
equivalence classes with respect to $\mathcal{D}$. Then
there is a bijection between maximal nilpotent subsemigroups
of $S$ and collections $(\rho_1,\dots,\rho_k)$ of linear
orders on the classes $\mathcal{M}_1$,\dots, $\mathcal{M}_k$
respectively.
\end{theorem}

\begin{proof}
Consider the partial order $\rho_T$ on $\mathcal{M}$, which
corresponds to some nilpotent subsemigroup
$T\subset S/R(S)$. If the idempotents $e,f\in\mathcal{M}$
belong to different $\mathcal{D}$-classes, then Lemma~\ref{lemma6}
implies that $exf=0$ for any $x\in S/R(S)$. Hence
$(e,f)\not\in\rho_T$. Therefore the elements from different
$\mathcal{M}_i$'s are not comparable with respect to $\rho_T$.
This means that $\rho_T$ decomposes into the disjoint union
$\rho_T=\rho_T^1\cup\dots\cup\rho_T^k$ of partial orders on the 
classes $\mathcal{M}_1$,\dots, $\mathcal{M}_k$
respectively.

Let now $(\rho_1^{(1)},\dots,\rho_k^{(1)})$ and
$(\rho_1^{(2)},\dots,\rho_k^{(2)})$ be two collections
of partial orders on the 
classes $\mathcal{M}_1$,\dots, $\mathcal{M}_k$. Set
$\rho_i=\cup_{j=1}^k\rho_j^{(i)}$, $i=1,2$. Let us show
that $\rho_1\subsetneqq\rho_2$ implies
$\mathrm{Mon}(\rho_1)\subsetneqq\mathrm{Mon}(\rho_2)$.
Indeed, that $\mathrm{Mon}(\rho_1)\subset\mathrm{Mon}(\rho_2)$
follows from Proposition~\ref{prop5}\eqref{prop5.1}. Choose
now some $\mathcal{M}_r$ and $e,f\in \mathcal{M}_r$ such that
$(e,f)\not\in \rho_1$ but $(e,f)\in \rho_2$. According to
\eqref{eqnil}, there exists $x\in S/R(S)$ such that
$exf\neq 0$. Then for the element $y=exf$ we have
$eyf=eexff=exf\neq 0$. On the other hand, from
\eqref{eqord} and Proposition~\ref{prop2} we obtain that
for any other pair $(g,h)\neq(e,f)$ of idempotents
from $\mathcal{M}$ we have $gyh=gexfh=0$
(indeed, if $g\neq e$ then $ge=0$, and if
$h\neq f$ then $fh=0$). Hence $y\in \mathrm{Mon}(\rho_2)$
but $y\not\in \mathrm{Mon}(\rho_1)$.

Thus the partial order $\rho_T$, which corresponds to
a maximal nilpotent subsemigroup $T\subset S/R(S)$ must
induce a maximal partial order (that is a linear order) 
on every class $\mathcal{M}_i$, $i=1,\dots,k$.

On the other hand, for an arbitrary collection,
$(\rho_1,\dots,\rho_k)$, of linear orders on the classes
$\mathcal{M}_1$,\dots, $\mathcal{M}_k$ respectively the
semigroup $\mathrm{Mon}(\rho)$, which is defined by the
partial order $\rho=\cup_{i=1}^k\rho_k$, is a maximal
nilpotent subsemigroup of $S$. Indeed,  
Proposition~\ref{prop5}\eqref{prop5.1} implies that
$\mathrm{Mon}(\rho)$ is nilpotent. If it is not maximal,
let $T\subset S/R(S)$ be any other nilpotent subsemigroup,
which properly contains $\mathrm{Mon}(\rho)$. Then the
partial order $\rho_T$ induces some partial
order $\rho_i^T$ on $\mathcal{M}_i$ for all $i=1,\dots,k$.
Moreover, from $\mathrm{Mon}(\rho)\subset T$ it follows
that $\rho_i\subset\rho_i^T$ for every $i$. Since all
$\rho_i$ are linear, we obtain $\rho_i=\rho_i^T$ for all $i$
and thus $\rho=\rho_T$. From Proposition~\ref{prop6}\eqref{prop6.1}
we now derive $\mathrm{Mon}(\rho)=\mathrm{Mon}(\rho_T)\supset T$,
a contradiction. 

Hence the maps $T\mapsto\rho_T$ and $\rho\mapsto\mathrm{Mon}(\rho)$
can be restricted to the sets $\mathrm{Nil}_{\mathrm{max}}(S)$ of all 
maximal nilpotent subsemigroups of $S$ and $\mathrm{Lin}$
of all collections of linear orders on the classes 
$\mathcal{M}_i$, $i=1,\dots,k$. Let us show that these restrictions
are, in fact, mutually inverse bijections.

If the collections $(\rho_1^{(1)},\dots,\rho_k^{(1)})$
and $(\rho_1^{(2)},\dots,\rho_k^{(2)})$ of linear orders
are different, without loss of generality we can assume that there
exist $e,f\in\mathcal{M}_1$ such that $(e,f)\in\rho_1^{(1)}$
but $(e,f)\not\in\rho_1^{(2)}$. Using the argument, analogous
to the one used in the proof of the implication
``$\rho_1\subsetneqq\rho_2$ implies
$\mathrm{Mon}(\rho_1)\subsetneqq\mathrm{Mon}(\rho_2)$'', one proves
that there exists $y\in S$ such that 
$y\in \mathrm{Mon}(\rho_1^{(1)}\cup\dots\cup\rho_k^{(1)})$
but $y\not\in \mathrm{Mon}(\rho_1^{(2)}\cup\dots\cup\rho_k^{(2)})$.
Hence the map $\rho\mapsto\mathrm{Mon}(\rho)$ is injective.

If $T_1$ and $T_2$ are two maximal nilpotent subsemigroups and
$\rho_{T_1}=\rho_{T_2}=\rho$, then from 
Proposition~\ref{prop6}\eqref{prop6.1} it follows that both
$T_1$ and $T_2$ are contained in the nilpotent subsemigroup
$\mathrm{Mon}(\rho)$. Hence $T_1=T_2$ and the map
$T\mapsto\rho_T$ is also injective. 

Since both, $\mathrm{Nil}_{\mathrm{max}}(S)$ and $\mathrm{Lin}$,
are finite sets, and both, $\rho\mapsto\mathrm{Mon}(\rho)$
and $T\mapsto\rho_T$, are injective maps, we obtain that both
these maps are, in fact, bijection. That they are inverse
to each other follows then from Proposition~\ref{prop6}\eqref{prop6.2}
and Proposition~\ref{prop6}\eqref{prop6.3}. 

This proves the statement of the theorem for the semigroup $S/R(S)$.
For the semigroup $S$ the statement now follows using 
Corollary~\ref{cor3}. This completes the proof.
\end{proof}

\begin{corollary}\label{cor6}
Let $S$ and $\mathcal{M}$ be as in Theorem~\ref{theorem4}. Then
the number of different maximal nilpotent subsemigroups in
$S$ equals $|\mathcal{M}_1|!\cdot|\mathcal{M}_2|!\cdots|\mathcal{M}_k|!$.
\end{corollary}

\begin{proof}
Follows from Theorem~\ref{theorem4} and the fact that the number of
linear orders on an $n$-element set equals $n!$.
\end{proof}

\begin{corollary}\label{cor6prime}
Let $S$ and $\mathcal{M}$ be as in Theorem~\ref{theorem4}. If
$|\mathcal{M}_1|=\dots=|\mathcal{M}_k|=1$, then $R(S)$ is the
unique maximal nilpotent subsemigroup of $S$.
\end{corollary}

\begin{proof}
Let $T$ be a maximal nilpotent subsemigroup of $S/R(S)$. Since
$\rho_T$ is a strict order, and all classes $\mathcal{M}_i$,
$i=1,\dots,k$, consist of $1$ element each, we have that
$e\rho_T f$ implies  $\mathcal{D}_e\neq\mathcal{D}_f$. 
However, in the proof of  Theorem~\ref{theorem4} it is shown that 
elements from different $\mathcal{D}$-classes are not comparable.
Hence $\rho_T=\varnothing$. From the condition \eqref{cond2}
it then follows that $T=\{0\}$. This proves the statement.
\end{proof}

\begin{corollary}\label{cor7}
Let $S$ be a finite inverse semigroup, and let $k$ be the number
of minimal $\mathcal{D}$-classes of $S$. Assume that
the minimal $\mathcal{D}$-classes contain $n_1,\dots,n_k$ idempotents
respectively. Then the number of maximal nilpotent subsemigroups in
$S$ equals $n_1!n_2!\cdots n_k!$.
\end{corollary}

\begin{proof}
Follows from Theorem~\ref{theorem3} and Corollary~\ref{cor6}. 
\end{proof}

\begin{corollary}\label{cor8}
Assume that $S$ is a finite semigroup. Assume further that
$S/R(S)$ contains a non-empty set, $\mathcal{M}$, of minimal
idempotents, which satisfies \eqref{cond1}-\eqref{cond3}. Then 
$R(S)$ coinsides with the intersection of all maximal
nilpotent subsemigroups of $S$.
\end{corollary}

\begin{proof}
Since $R(S)$ is contained in each maximal nilpotent subsemigroup
of $S$ by Corollary~\ref{cor3}, it is enough to show that
the intersection of all maximal nilpotent subsemigroups of 
$S/R(S)$ is $\{0\}$. The latter follows from
Theorem~\ref{theorem4} and the obvious fact that for any set
$K$ the intersection of all linear orders on $K$ is the
empty relation.
\end{proof}

\begin{problem}
Describe all finite semigroups, the radical of which coincides with
the intersection of all maximal nilpotent subsemigroups. 
\end{problem}

For semigroups, whose radical is zero, Theorem~\ref{theorem4} can
be substantially streng\-then\-ed. For $m\in\mathbb{N}$ let 
$\mathrm{Nil}_m(S)$ denote the set of all nilpotent 
subsemigroups of $S$ of nilpotency class $\leq m$, and let 
$\mathrm{Ord}_m(\mathcal{M})$ denote the set of all (strict)
partial orders on $\mathcal{M}$ in which the cardinalities of
chains do not exceed $m$.

Let $M$ be a set. By an {\em ordered partition of $M$ into 
$l$ blocks} we will mean a partition, $M=M_1\cup\dots\cup
M_l$, into $l$ non-empty blocks in which the order of blocks 
is also taken into account. Each usual partition of $M$ 
into $l$ blocks gives, obviously, $l!$ ordered partitions. With every 
ordered partition $M=M_1\cup\dots\cup 
M_l$ we associate the set 
\begin{displaymath}
\mathrm{Ord}(M_1,\dots,M_l)=
\bigcup_{1\leq i<j\leq k} M_i\times M_j \subset
M\times M.
\end{displaymath}

\begin{lemma}\label{l33}{\rm (\cite[Lemma~7]{GM1})}
Fix a positive integer, $l\leq |M|$. Then for every ordered 
partition $M=M_1\cup\dots\cup M_l$ the 
set $\mathrm{Ord}(M_1,\dots,M_l)$ is a maximal
element in $\mathrm{Ord}_l(M)$. Different ordered partitions 
of $M$ correspond to different elements in 
$\mathrm{Ord}_l(M)$, and each maximal element in
$\mathrm{Ord}_l(M)$ has the form
$\mathrm{Ord}(M_1,\dots,M_l)$ for some ordered
partition $M=M_1\cup\dots\cup M_l$.
\end{lemma}

Consider now the set $\mathcal{M}$ with the fixed decomposition 
$\mathcal{M}=\mathcal{M}_1\cup\dots\cup\mathcal{M}_k$ into equivalence
classes with respect to the restriction of the $\mathcal{D}$-relation. 
Let $n=\mathrm{max}_i|\mathcal{M}_i|$ and $m\leq n$ be fixed. On every
class $\mathcal{M}_i$, $i=1,\dots,k$, we choose some partial order,
say $\rho_i$, according to the following rule: If $|\mathcal{M}_i|\geq m$,
we choose some ordered decomposition,
$\mathcal{M}_i=\mathcal{M}_i^{(1)}\cup\dots\cup\mathcal{M}_i^{(m)}$, of
$\mathcal{M}_i$ into $m$ blocks, and define $\rho_i$ as
$\mathrm{Ord}(\mathcal{M}_i^{(1)},\dots,\mathcal{M}_i^{(m)})$. 
If $|\mathcal{M}_i|< m$, we choose
some linear order on $\mathcal{M}_i$. Set $\rho=\rho_1\cup\dots\cup\rho_k$
and denote by $\mathrm{Ord}^{\mathrm{part}}_m(\mathcal{M})$ the
set of all partial orders on $\mathcal{M}$, which can be obtained in this
way.

\begin{theorem}\label{theorem5}
Let $S$ be a finite semigroup such that $R(S)=\{0\}$, 
$\mathcal{M}$ be a set of minimal idempotents of $S$, satisfying
\eqref{cond1}-\eqref{cond3}, 
$\mathcal{M}=\mathcal{M}_1\cup\dots\cup\mathcal{M}_k$ be the decomposition
of $\mathcal{M}$ into equivalence
classes with respect to the restriction of the $\mathcal{D}$-relation,
and $n=\mathrm{max}_i|\mathcal{M}_i|$. Then the nilpotency class of every
nilpotent semigroup in $S$ does not exceed $n$, the nilpotency class of 
every maximal nilpotent semigroup in $S$ equals $n$, and for every
$m\leq n$ the maps $T\mapsto\rho_T$ and $\rho\mapsto \mathrm{Mon}(\rho)$ 
induce mutually inverse bijections between the set of all maximal
nilpotent subsemigroups of $S$ of nilpotency class $\leq m$
(i.e. the set of all maximal elements in $\mathrm{Nil}_m(S)$) 
and $\mathrm{Ord}^{\mathrm{part}}_m(\mathcal{M})$.
\end{theorem}

\begin{proof}
First we claim that $T\in \mathrm{Nil}_m(S)$ implies
$\rho_T\in \mathrm{Ord}_m(\mathcal{M})$. Indeed,
let $e_1\rho_1e_2\rho_2e_3\cdots\rho_Te_p$ be a chain of cardinality
$p$ in $\rho_T$. Then there exist $x_1,\dots,x_{p-1}\in T$ such that
$e_1x_1e_2\neq 0$,\dots, $e_{p-1}x_{p-1}e_p\neq 0$. Applying
\cite[Proposition~2.4]{Gri} and \eqref{cond3} in the same way as it
was done in the proof of Proposition~\ref{prop5}\eqref{prop5.2},
we obtain that $e_1x_1x_2\cdots x_{p-1}e_p\neq 0$ and hence
$x_1x_2\cdots x_{p-1}\neq 0$. Thus $p-1\leq m-1$ and so $p\leq m$.
This means that $\rho_T\in \mathrm{Ord}_m(\mathcal{M})$.

The arguments, analogous to those used in the proof of
Proposition~\ref{prop5}\eqref{prop5.1}, show that
$\rho\in \mathrm{Ord}_m(\mathcal{M})$ implies
$\mathrm{Mon}(\rho)\in \mathrm{Nil}_m(S)$. 

Now, analogously to the proof of Theorem~\ref{theorem4} one
shows that the map $T\mapsto \rho_T$ maps the set
$\mathrm{Nil}_m^{\mathrm{max}}(S)$ of the maximal elements in
$\mathrm{Nil}_m(S)$ to the set 
$\mathrm{Ord}^{\mathrm{part}}_m(\mathcal{M})$,  the map
$\rho\mapsto \mathrm{Mon}(\rho)$ maps 
$\mathrm{Ord}^{\mathrm{part}}_m(\mathcal{M})$ to
$\mathrm{Nil}_m^{\mathrm{max}}(S)$, and that the restrictions
of these maps to these sets are injective. In particular,
$|\mathrm{Nil}_m^{\mathrm{max}}(S)|=
|\mathrm{Ord}^{\mathrm{part}}_m(\mathcal{M})|$. It further
follows from Proposition~\ref{prop6}\eqref{prop6.2}
and Proposition~\ref{prop6}\eqref{prop6.3} that the above maps
induce mutually inverse bijections between 
$\mathrm{Nil}_m^{\mathrm{max}}(S)$ and
$\mathrm{Ord}^{\mathrm{part}}_m(\mathcal{M})$.

Finally, from Theorem~\ref{theorem4} it follows that the
maximal nilpotent subsemigroups of $S$ correspond to the
elements of $\mathrm{Ord}^{\mathrm{part}}_n(\mathcal{M})$. Hence
their nilpotency class equals $n$, and the nilpotency class of
any nilpotent subsemigroup does not exceed $n$. This completes
the proof.
\end{proof}

\begin{corollary}\label{cor9}
Let $S$ be as in Theorem~\ref{theorem5} and $m\leq n$.
\begin{enumerate}[(a)]
\item\label{cor9.1}
$S$ contains exactly
\begin{displaymath}
\prod_{i:|\mathcal{M}_i|\geq m}
\left(\sum_{j=0}^{m-1}(-1)^j\binom{m}{j}(m-j)^{|\mathcal{M}_i|}\right)
\cdot \prod_{i:|\mathcal{M}_i|< m}|\mathcal{M}_i|!
\end{displaymath}
maximal nilpotent subsemigroups of nilpotency class $m$.
\item\label{cor9.2} If a maximal nilpotent subsemigroup, 
$T\subset S$, of nilpotency  class $m$,
is determined by the ordered partitions
$\mathcal{M}_i=\mathcal{M}_i^{(1)}\cup\dots\cup\mathcal{M}_i^{(m)}$
of the sets $\mathcal{M}_1$,\dots, $\mathcal{M}_k$, then
$T$ is contained in exactly 
\begin{displaymath}
\prod_{i=1}^k\prod_{j=1}^m |\mathcal{M}_i^{(j)}|!
\end{displaymath}
maximal nilpotent subsemigroups of $S$.
\end{enumerate}
\end{corollary}

\begin{proof}
By Theorem~\ref{theorem5} the number of 
maximal nilpotent subsemigroups of $S$ of nilpotency class $m$
equals $|\mathrm{Ord}^{\mathrm{part}}_m(\mathcal{M})|$. Each partial
order $\rho\in \mathrm{Ord}^{\mathrm{part}}_m(\mathcal{M})$ is
determined by ordered decompositions into $m$ blocks of those
$\mathcal{M}_i$'s, which contain at least $m$ elements
(for such $\mathcal{M}_i$ the number of ordered partitions into
$m$ blocks equals the number of surjections from $\mathcal{M}_i$
to $\{1,\dots,m\}$ and hence equals 
$\sum_{j=0}^{m-1}(-1)^j\binom{m}{j}(m-j)^{|\mathcal{M}_i|}$),
and linear orders on those $\mathcal{M}_i$'s, which contain less
than $m$ elements (such $\mathcal{M}_i$ can be linearly ordered in
$|\mathcal{M}_i|!$ different ways). The ordered partitions of
$\mathcal{M}_i$'s of the first type and the linear orders
on $\mathcal{M}_i$'s of the second type can be chosen independently,
which proves \eqref{cor9.1}.

To prove \eqref{cor9.2} we note that from Theorem~\ref{theorem5}
and Proposition~\ref{prop5} it follows that the number of those 
maximal nilpotent subsemigroups of $S$, which contain $T$, equals 
the number of extensions of the partial order
$\rho_T\in\mathrm{Ord}^{\mathrm{part}}_m(\mathcal{M})$ (which is
defined by the ordered partitions
$\mathcal{M}_i=\mathcal{M}_i^{(1)}\cup\dots\cup\mathcal{M}_i^{(m)}$
of the sets $\mathcal{M}_1$,\dots, $\mathcal{M}_k$) to a
partial order from $\mathrm{Ord}^{\mathrm{part}}_n(\mathcal{M})$.
It is clear that for such an extension one has to choose
(independently) linear orders on every $\mathcal{M}_i^{(j)}$.
This can be done in $\prod_{i=1}^k\prod_{j=1}^m |\mathcal{M}_i^{(j)}|!$
different ways, and the proof is complete.
\end{proof}

\section{Examples and applications}\label{s6}

In this section we show that in almost all cases, where 
maximal nilpotent semigroups are classified in the literature,
this can be done using Theorem~\ref{theorem4} or
Theorem~\ref{theorem5}. The only exception which we know is the
semigroup of all linear operators on a finite-dimensional
vector space over a finite field. The easiest argument is
for finite inverse semigroups. As we have already seen, the
radical of such semigroups is always $\{0\}$, and the
conditions \eqref{cond1}-\eqref{cond3} are always satisfied
for the set of all minimal idempotents. Hence to such semigroups
Theorem~\ref{theorem5} and Corollary~\ref{cor9} can be applied
immediately. In the sequel for a positive integer, $n$, we set
$N=\{1,2,\dots,n\}$ and we also denote by $S_n$ the symmetric 
group on $N$. For transformation semigroups we will use
the {\em right action} notation and will denote by $x\,\alpha$
the image of $x$ under $\alpha$. By a {\em proper}
subset of a set, $X$, we will mean a subset of $X$, different
from $X$ and $\varnothing$.

\subsection{The symmetric inverse semigroup}\label{s6.1}

The symmetric inverse semigroup $\mathcal{IS}_n$, which consists
of all partial injective transformations of the set $N$,
contains a unique minimal $\mathcal{D}$-class, namely $D_1$, which
consists of all transformations of rank $1$. $D_1$ contains
$n$ idempotents, which are the identity transformations on the
one-element subsets of $\{1,2,\dots,n\}$. Hence $\mathcal{IS}_n$ has
$n!$ maximal nilpotent subsemigroups (each of nilpotency class $n$)
and $\sum_{i=0}^{k-1}(-1)^i\binom{k}{i}(k-i)^n$ maximal nilpotent
subsemigroups of a given nilpotency class, $k$. This is proved
in \cite{GK1,GK2}.

\subsection{The semigroup $\mathcal{IO}_n$}\label{s6.105}

The semigroup $\mathcal{IO}_n$ of all injective order-preserving 
transformations of the set $N$ is an inverse semigroup and it
also contains a unique minimal $\mathcal{D}$-class, which coincides 
with the  class $D_1$ of the semigroup $\mathcal{IS}_n$. Hence
$\mathcal{IO}_n$ contains  $n!$ maximal nilpotent subsemigroups 
(each of nilpotency class $n$) and 
$\sum_{i=0}^{k-1}(-1)^i\binom{k}{i}(k-i)^n$ maximal nilpotent
subsemigroups of a given nilpotency class, $k$. This is proved
in \cite[Theorem~7]{GM1}.

\subsection{The semigroup $\mathrm{PAut}(\mathbb{F}^n_q)$}\label{s6.2}

The semigroup $\mathrm{PAut}(\mathbb{F}^n_q)$ of all partial 
automorphisms of the $n$-dimensional vector space
$\mathbb{F}^n_q$ over the field $\mathbb{F}_q$ with $q<\infty$
elements is a finite inverse semigroup, which
consists of all linear isomorphisms $\varphi:U\to V$,
where $U$ and $V$ are arbitrary subspaces of $\mathbb{F}^n_q$
(of the same dimension). $\mathrm{PAut}(\mathbb{F}^n_q)$
has a unique minimal $\mathcal{D}$-class, namely $D_1$, 
which consists of all elements of rank $1$, that is 
linear isomorphisms between subspaces of dimension $1$. The
idempotents in $D_1$ are the identity transformations of the 
one-dimensional subspaces of $\mathbb{F}^n_q$. Since
$\mathbb{F}^n_q$ contains $m=\frac{q^n-1}{q-1}$ one-dimensional
subspaces, we obtain that $\mathrm{PAut}(\mathbb{F}^n_q)$
has $m!$ maximal nilpotent subsemigroups (each of nilpotency class 
$m$). This was proved in \cite{Sh1,Sh2}. Moreover, for every 
$k\leq m$, $\mathrm{PAut}(\mathbb{F}^n_q)$
contains $\sum_{i=0}^{k-1}(-1)^i\binom{k}{i}(k-i)^m$ maximal 
nilpotent subsemigroups of nilpotency class $k$. 

\subsection{The semigroup $B_n$ of binary relations}\label{s6.3}

The semigroup $B_n$ of all binary relations on the set
$N$ is not inverse. The zero element of this
semigroup is the empty relation. Let $\theta$ denote the
full relation on $N$. Then for every non-empty relation
$\alpha\in B_n$ one obviously has $\theta\alpha\theta=\theta$,
which implies that every non-zero ideal of $B_n$ contains
the idempotent $\theta$ and hence is not nilpoptent. In particular,
$R(B_n)=\{0\}$. Furthermore, from the above argument
(or, alternatively, using \cite{Za,PW}) one also obtains that
$B_n$ contains a unique minimal $\mathcal{D}$-class, namely
$D=\mathcal{D}_{\theta}$. From  \cite{Za} we have $|D|=(2^n-1)^2$
and $\varphi\in D$ if and only if $\varphi=A\times B$, where
$A$ and $B$ are two non-empty subsets of $N$. The class
$D$ contains idempotents $\varepsilon_i=\{(i,i)\}$, $i\in N$,
which satisfy the conditions \eqref{cond1}-\eqref{cond3}
(for \eqref{cond1} and \eqref{cond2} this is obvious, and 
\eqref{cond3} follows from the fact that $\varepsilon_i\alpha
\varepsilon_j$ is either $\varnothing$ or $\{(i,j)\}$, moreover,
if $\varepsilon_i\alpha\beta\varepsilon_j=\{(i,j)\}$ then 
$(i,k)\in\alpha$ and $(k,j)\in\beta$ for some $k$ and hence
$\varepsilon_i\alpha\varepsilon_k\beta\varepsilon_j=\{(i,j)\}$).

\begin{theorem}\label{theorem7}
\begin{enumerate}[(a)]
\item\label{theorem7.1}
The semigroup $B_n$ contains $n!$ maximal nilpotent subsemigroups,
each of nilpotency class $n$.
\item\label{theorem7.2} For every $k\in\mathbb{N}$,
$k\leq n$, $B_n$ contains $\sum_{i=0}^{k-1}(-1)^i\binom{k}{i}(k-i)^n$
maximal nilpotent subsemigroups of nilpotency class $k$.
\item\label{theorem7.3} If $T_1$ and $T_2$ are two 
maximal nilpotent subsemigroups of $B_n$ then there exists
a permutation, $\pi\in S_n$, such that $T_2=\pi^{-1}T_1\pi$, in
particular, all maximal nilpotent subsemigroups of $B_n$ are
isomorphic.
\item\label{theorem7.4} Every maximal nilpotent subsemigroup of
$B_n$ consists of $2^{n(n-1)/2}$ elements.
\end{enumerate}
\end{theorem}

\begin{proof}
The statement \eqref{theorem7.1} follows from Theorem~\ref{theorem5}
applied to the set $\{\varepsilon_1,\dots,\varepsilon_n\}$ of
minimal idempotents. The statement \eqref{theorem7.2} follows
from Corollary~\ref{cor9}. It is obvious that the transformation of
a given linear order on $\{\varepsilon_1,\dots,\varepsilon_n\}$ into
another linear order is determined by some permutation, $\pi$, 
of the elements of $N$. Let $T_1$ be the maximal nilpotent subsemigroup,
which corresponds to the original linear order, and $T_2$ be the
maximal nilpotent subsemigroup, which corresponds to the new linear order.
Then for any $\alpha\in B_n$ we have
\begin{displaymath}
\varepsilon_i\alpha\varepsilon_j\neq 0\quad\Leftrightarrow\quad
\varepsilon_{\pi(i)}\cdot\pi^{-1}\alpha\pi\cdot\varepsilon_{\pi(j)}
\neq 0,
\end{displaymath}
which implies $T_2=\pi^{-1}T_1\pi$. Since $\alpha\mapsto \pi^{-1}\alpha\pi$
is an automorphism of $B_n$, we derive that $T_1$ and $T_2$ are isomorphic,
which gives \eqref{theorem7.3}.

Finally, from \eqref{theorem7.3} it follows that to prove \eqref{theorem7.4}
we have to compute the cardinality of the maximal nilpotent subsemigroup
$T$ of $B_n$, which corresponds to the natural order 
$\varepsilon_1<\varepsilon_2<\dots<\varepsilon_n$. The semigroup $T$
consists of all $\alpha\in B_n$ such that 
$\varepsilon_i\alpha\varepsilon_j\neq 0$ implies $i<j$. However,
$\varepsilon_i\alpha\varepsilon_j\neq 0$ if and only if 
$(i,j)\in\alpha$. Hence the elements of $T$ are just arbitrary
subsets of the set $\{(i,j):1\leq i<j\leq n\}$. Therefore
$|T|=2^{n(n-1)/2}$ and the proof is complete.
\end{proof}

\begin{remark}\label{rem100}
{\rm
We note that, since $B_n$ does not act naturally on the set $N$, one can 
not apply to $B_n$ (at least in a straightforward way) the general 
approach to the description of  maximal nilpotent subsemigroups 
of transformation semigroups, developed in \cite[Section~6]{GM1}.
}
\end{remark}

\begin{remark}\label{rem101}
{\rm
We also note that in \cite[Corollary~2]{GM3} it was shown that the 
semigroup  $B_n$ is almost nilpotent in the sense that it contains a 
nilpotent subsemigroup, $R_n$, such that $\frac{|R_n|}{|B_n|}\to 1$,
$n\to\infty$. However, the zero in $R_n$ is the full relation
$\theta$ and not the empty relation $0$.
}
\end{remark}

\subsection{The semigroup $D_n$}\label{s6.5}

The semigroup $D_n$ consists of all 
order-decreasing transformations of $N$, that
is of all transformations $\alpha$ such that $x\,\alpha\leq x$
for all $x\in N$. The zero element of this semigroup is
the transformation $0:x\mapsto 1$ for all $x\in N$. A
transformation, $\alpha\in D_n$ is nilpotent if and only if
$x\, \alpha<x$ for all $x>1$. It is obvious that for a nilpotent
element, $\alpha$, the ideal $D_n^1\alpha D_n^1$
contains only nilpotent elements and hence, by Proposition~\ref{prop1},
is nilpotent. This means that the radical
$R({D}_n)$ coincides with the set of all nilpotent
elements (and hence with the unique maximal nilpotent subsemigroup 
of ${D}_n$). It also follows that the Rees factor
$S={D}_n/R({D}_n)$ contains only the trivial
nilpotent subsemigroup $\{0\}$.

Let us also show that ${D}_n$ satisfies the conditions of
Corollary~\ref{cor6prime}. It is obvious that 
$\alpha\in {D}_n/R({D}_n)$ if and only if there
exists $k>1$ such that $k\,\alpha=k$. This means that every 
element of rank two in $S$ is an idempotent. It is also easy to see
that the minimal $\mathcal{D}$-classes of $S$ consist of 
the elements of rank two and that the idempotents $e$ and $f$
are $\mathcal{D}$-related if and only if 
$\mathrm{im}(e)=\mathrm{im}(f)$.

Let now $e,f,g$ be idempotents of $S$ such that 
$\mathrm{im}(e)=\mathrm{im}(f)=\{1,k\}$ and $\mathrm{im}(g)=\{1,m\}$,
where $k\neq m$. Then $ef=e$ and $eg\in R({D}_n)$. Hence
two minimal idempotents of $S$ commute if and only if they belong
to different $\mathcal{D}$-classes. Let us now for $k=2,\dots,n$,
choose some idempotent, $e_k$, in the $\mathcal{D}$-class, defined
by the image set $\{1,k\}$. Then the set $\mathcal{M}=
\{e_2,\dots,e_n\}$ satisfies \eqref{cond1}. The condition
\eqref{cond2} is also satisfied since for $0\neq\alpha\in S$
and $k>1$ we have that $k\,\alpha=k$ implies
$e_k\alpha e_k=e_k\neq 0$. 

Finally, assume that $e_k\alpha e_m\beta e_l\neq 0$. Then 
there exists $x>1$ such that $x\,e_k\alpha e_m\beta e_l=x$.
Since all elements in $S$ are order-decreasing, the latter
equality is possible if and only if 
$x=k=k\,\alpha=m=m\,\beta=l=x$. This implies that
$e_k\alpha e_m\beta e_l=e_k\alpha\beta e_l$ and the condition
\eqref{cond3} follows.

Let us note that the nilpotent subsemigroups of $D_n$ were
deeply studied in \cite{St}. In particular, it is shown that
all maximal nilpotent subsemigroups of $D_n$ of a fixed
nilpotency class are not isomorphic.

\subsection{The semigroup $\mathcal{PT}_n$}\label{s6.6}

Let us now consider the semigroup $\mathcal{PT}_n$ of all
partial transformations of $N$. This semigroup
has a unique minimal $\mathcal{D}$-class, namely the one which
consists of all transformations of rank $1$. Let $\varepsilon_i$,
$i=1,\dots,n$ denote the minimal idempotent, defined via
$\mathrm{dom}(\varepsilon_i)=\mathrm{im}(\varepsilon_i)=\{i\}$.
In the same way as for the semigroup $B_n$ (see Subsection~\ref{s6.3})
one shows that the set $\mathcal{M}=\{1,2,\dots,n\}$ satisfies
the conditions \eqref{cond1}-\eqref{cond3}.

\begin{theorem}\label{theorem8}
\begin{enumerate}[(a)]
\item\label{theorem8.1}
The semigroup $\mathcal{PT}_n$ contains $n!$ maximal nilpotent
subsemigroups, each of nilpotency class $n$.
\item\label{theorem8.2}
For every positive integer $k<n$, the semigroup
$\mathcal{PT}_n$ contains exactly 
$\sum_{i=0}^{k-1}(-1)^i\binom{n}{i}(k-i)^n$
maximal nilpotent subsemigroups of nilpotency class $k$.
\item\label{theorem8.3} If $T_1$ and $T_2$ are two maximal 
nilpotent subsemigroups of $\mathcal{PT}_n$ then there exists
a permutation, $\pi\in S_n$, such that $T_2=\pi^{-1}T_1\pi$, in
particular, all maximal nilpotent subsemigroups of $B_n$ are
isomorphic.
\item\label{theorem8.4} Every maximal nilpotent subsemigroup of
$\mathcal{PT}_n$ consists of $n!$ elements.
\end{enumerate}
\end{theorem}

\begin{proof}
The statements \eqref{theorem8.1}, \eqref{theorem8.2}, and
\eqref{theorem8.3} are proved mutatis mutandis the correspondent
statements of Theorem~\ref{theorem7}. To prove \eqref{theorem8.4}
we observe that the condition $\varepsilon_i\alpha\varepsilon_j\neq 0$
is equivalent to the fact that $i\in\mathrm{dom}(\alpha)$ and
$i\,\alpha=j$. Consider now the linear order
$\varepsilon_1<\varepsilon_2<\dots<\varepsilon_n$. All implications
\begin{displaymath}
\varepsilon_i\alpha\varepsilon_j\neq 0\quad\Rightarrow\quad
i<j
\end{displaymath}
are true if and only if for every $i\in\mathrm{dom}(\alpha)$
one has $i\,\alpha>i$. Since the values of $\alpha$ at different
elements of $N$ can be chosen independently, we get exactly 
$n!$ different possibilities for $\alpha$. This completes the proof.
\end{proof}

\subsection{The symmetric semigroup $\mathcal{T}(X)$}\label{s6.7}

Let $X$ be a set. Then the semigroup $\mathcal{T}(X)$ of all 
transformations of $X$ does not contain the zero element. 
However, one can consider the nilpotent subsemigroups of 
$\mathcal{T}(X)$, whose zero coincides with a fixed idempotent, 
$e\in \mathcal{T}(X)$. From this point of view nilpotent 
subsemigroups of $\mathcal{T}(X)$ were studied in \cite{BRR}.

Let $|X|=n$ and $e$ be a fixed idempotent of $\mathcal{T}(X)$,
$\mathrm{im}(e)=\{a_1,\dots,a_k\}$, 
$A_i=\{x\in X: x\, e=a_i\}$, $i=1,\dots,k$.

\begin{lemma}\label{lemma10}
The idempotent $e$ is a two-sided zero for an element, 
$\alpha\in \mathcal{T}(X)$, if and only if $a_i\,\alpha=a_i$
and $A_i\,\alpha\subset A_i$ for all $i=1,\dots,k$.
\end{lemma}

\begin{proof}
A direct computation.
\end{proof}

Let $A$ be an arbitrary non-empty set and $a$ be a 
fixed element of $A$. Then the map
\begin{displaymath}
\begin{array}{ccc}
S_a=\{\alpha\in\mathcal{T}(A):a\,\alpha=a\} & \rightarrow &
\mathcal{PT}(A\setminus\{a\})\\
\alpha & \mapsto & \tilde{\alpha},
\end{array}
\end{displaymath}
where
\begin{displaymath}
x\,\tilde{\alpha}=
\begin{cases}
x\,\alpha, & x\,\alpha\neq a\\
\text{not defined}, & x\,\alpha=a, 
\end{cases}
\end{displaymath}
is an isomorphism. From this and Lemma~\ref{lemma10} we obtain
that the semigroup
$\mathcal{T}_e(X)$ of all elements from $\mathcal{T}(X)$, for which
$e$ is the two-sided zero, is isomorphic to the direct product
$\mathcal{PT}(A_1\setminus\{a_1\})\times
\dots\times\mathcal{PT}(A_k\setminus\{a_k\})$.

Since all factors of the latter product are semigroups with zero,
each maximal nilpotent subsemigroup of the product decomposes into 
a direct product of some maximal nilpotent subsemigroups of
the factors. From the latter argument and Theorem~\ref{theorem8} 
(or, alternatively, from the fact that the set of all minimal
$\mathcal{D}$-classes  of a direct product can be identified with
the union of all minimal $\mathcal{D}$-classes of factors; and
a subsequent application of Theorem~\ref{theorem5}) we obtain:

\begin{theorem}{\rm
(\cite[Theorem~4.1 and Corollary~4.3]{BRR})}\label{theorem9}
The semigroup $\mathcal{T}(X)$ contains 
$|A_1\setminus\{a_1\}|!\cdot\cdots\cdot|A_k\setminus\{a_k\}|!$
maximal nilpotent subsemigroups for which $e$ is the zero element. 
These semigroups
are in bijection with the collections of linear orders on the sets
$A_1\setminus\{a_1\}$,\dots, $A_k\setminus\{a_k\}$. They all are
pairwise isomorphic and contain
$|A_1\setminus\{a_1\}|!\cdots|A_k\setminus\{a_k\}|!$
elements each.
\end{theorem}

Let us note that Theorem~\ref{theorem8} even tells us 
that an isomorphism between two maximal nilpotent subsemigroups from
$\mathcal{T}(X)$ for which $e$ is the zero element can be constructed 
as a conjugation by some permutation on $X$, which acts identically on 
$\mathrm{im}(e)$.

\subsection{Variants of $\mathcal{IS}_n$}\label{s6.8}

For a fixed element, $a$, of a semigroup, $S$, one defines
a new associative {\em sandwich operation} $*_a$ on $S$ via
$x*_ay=xay$. If $a$ is fixed, we will simply use the notation
$x*y$. The semigroup $(S,{*_a})$ is called a 
{\em variant} of $S$.

The nilpotent subsemigroups of the variants of the
symmetric inverse semigroup $\mathcal{IS}_n$ were studied in
\cite{Ts}. In particular, the structure of maximal nilpotent 
subsemigroups of a given nilpotency class, $k$, was
described in \cite[Theorem~14]{Ts}. From this description
we immediately get the following:

\begin{proposition}\label{prop10}
Let $\varepsilon$ be an element of $\mathcal{IS}_n$ of rank $k$.
Then the variant $(\mathcal{IS}_n,*_{\varepsilon})$
contains $k!$ maximal nilpotent subsemigroups.
\end{proposition}

Let us show that this can be obtained from Theorem~\ref{theorem4}.
In \cite{Ts2} it is shown that 
$(\mathcal{IS}_n,*_{\alpha})\cong
(\mathcal{IS}_n,*_{\beta})$ if and only if 
$\mathrm{rank}(\alpha)=\mathrm{rank}(\beta)$, which allows us
to assume that $\varepsilon$ is an idempotent. Further we
can assume that $k<n$ since for $k=n$ we get
$(\mathcal{IS}_n,*_{1})=\mathcal{IS}_n$.
Set $\mathcal{I}_{\varepsilon}:=(\mathcal{IS}_n,*_{\varepsilon})$.
The following statement is obvious:

\begin{lemma}\label{l11}
The element $\alpha\in \mathcal{IS}_n$ is an idempotent in
$\mathcal{I}_{\varepsilon}$ if and only if $\alpha$ is
an idempotent of $\mathcal{IS}_n$ and $\mathrm{dom}(\alpha)\subset
\mathrm{dom}(\varepsilon)$.
\end{lemma}

\begin{corollary}\label{cor11}
\begin{enumerate}[(a)]
\item\label{cor11.1}
Any two idempotents in $\mathcal{I}_{\varepsilon}$ commute.
\item\label{cor11.2}
An idempotent, $\alpha\in \mathcal{I}_{\varepsilon}$,
is primitive if and only if $\mathrm{rank}(\alpha)=1$.
\end{enumerate}
\end{corollary}

For $m\in\mathrm{dom}(\varepsilon)$ we denote by
$\varepsilon_m$ the primitive idempotent of $\mathcal{I}_{\varepsilon}$
satisfying $\mathrm{dom}(\varepsilon_m)=\{m\}$. According
to \cite[Theorem~2.2]{Ts3}, all minimal idempotents of
$\mathcal{I}_{\varepsilon}$ belong to the same
$\mathcal{D}$-class. From Theorem~\ref{theorem1} it then follows that
the minimal and the primitive idempotents in
$\mathcal{I}_{\varepsilon}$ coincide. In particular,
the set $\mathcal{M}=\{\varepsilon_m:m\in\mathrm{dom}(\varepsilon)\}$
of all minimal idempotents of $\mathcal{I}_{\varepsilon}$
satisfies \eqref{cond1}.

\begin{lemma}\label{lemma11}
We have $R(\mathcal{I}_{\varepsilon})=\{\alpha\in
\mathcal{IS}_n:\varepsilon\alpha\varepsilon=0\}$.
\end{lemma}

\begin{proof}
If $\varepsilon\alpha\varepsilon=0$, then the definition of
$*_{\varepsilon}$ immediately implies that
\begin{displaymath}
(\alpha)=\mathcal{I}_{\varepsilon}^1*\alpha*\mathcal{I}_{\varepsilon}^1
= \mathcal{I}_{\varepsilon}^1\varepsilon\alpha\varepsilon
\mathcal{I}_{\varepsilon}^1=0,
\end{displaymath}
and hence $\alpha\in R(\mathcal{I}_{\varepsilon})$.

On the other hand, if $\varepsilon\alpha\varepsilon\neq 0$,
then for the element $\alpha^{-1}\in \mathcal{IS}_n$ we have
\begin{displaymath}
\varepsilon*\alpha*\alpha^{-1}*\varepsilon=
\varepsilon^2\alpha\varepsilon\alpha^{-1}\varepsilon^2=
(\varepsilon\alpha\varepsilon)\cdot(\varepsilon\alpha^{-1}\varepsilon)=
(\varepsilon\alpha\varepsilon)\cdot
(\varepsilon\alpha\varepsilon)^{-1}=\gamma.
\end{displaymath}
The element $\gamma$ is a non-zero idempotent in
$\mathcal{IS}_n$, whose domain is contained in $\mathrm{dom}(\varepsilon)$.
By Lemma~\ref{lemma11}, $\gamma$ a non-zero idempotent of
$\mathcal{I}_{\varepsilon}$ as well. As $\gamma\in (a)$, we obtain that
$\alpha\not\in R(\mathcal{I}_{\varepsilon})$ and the statement is proved.
\end{proof}

If $\varepsilon\alpha\varepsilon\neq 0$, $m\in \mathrm{dom}
(\varepsilon\alpha\varepsilon)$ and $l=m\,\alpha$, we have
$\varepsilon_m*\alpha*\varepsilon_l\neq 0$. Hence the set 
$\mathcal{M}$ satisfies \eqref{cond2}.

Finally, $\varepsilon_m*\alpha*\beta*\varepsilon_l\neq 0$ if and only
if $m\in \mathrm{dom}(\alpha)$, $m\, \alpha\in \mathrm{dom}(\varepsilon)
\cap \mathrm{dom}(\beta)$ and $l=m\,\alpha\beta$. On the other hand,
this implies that $\varepsilon_m*\alpha*\varepsilon_p*\beta*\varepsilon_l
\neq 0$ if and only if we have that the same conditions are satisfied
and $p=x\,\alpha$. Moreover, in this case
\begin{displaymath}
\varepsilon_m*\alpha*\beta*\varepsilon_l=
\varepsilon_m*\alpha*\varepsilon_p*\beta*\varepsilon_l.
\end{displaymath}
Hence the condition \eqref{cond3} for $\mathcal{M}$ is satisfied as well.
Since all elements from $\mathcal{M}$ belong to the same
$\mathcal{D}$-class, we obtain that Proposition~\ref{prop10}
follows from Theorem~\ref{theorem4}.

\subsection{The factor power $\mathcal{FP}^+(S_n)$}\label{s6.9}

Let $\mathcal{P}(S_n)$ be the power semigroup of $S_n$, that is
the set $\{A:A\subset S_n\}$ with the natural operation
$A\cdot B=\{\alpha\beta:\alpha\in A,\beta\in B\}$. On the 
semigroup $\mathcal{P}(S_n)$ we define the equivalence relation
$\sim$ as follows: for $A,B\in \mathcal{P}(S_n)$ we have
$A\sim B$ if and only if for each $i\in N$ the sets
$\{i\,\alpha:\alpha\in A\}$ and $\{i\,\beta:\beta\in B\}$
coincide. In is straightforward to verify that $\sim$
is a well-defined congruence on $\mathcal{P}(S_n)$. The
corresponding quotient semigroup
$\mathcal{P}(S_n)/\sim$ is called the {\em factor power}
of $S_n$ and is denoted by $\mathcal{FP}(S_n)$. The semigroup
$\mathcal{FP}(S_n)$ has an adjoint zero element, which is
the class, consisting of the empty subset of $S_n$. 
Taking this class away, we obtain the semigroup
$\mathcal{FP}^+(S_n)$, which we will also call the
{\em factor power} of $S_n$, abusing the language. In this section
we will consider the semigroup $\mathcal{FP}^+(S_n)$ only.

If $A\subset S_n$, then with the corresponding class
$\overline{A}\in \mathcal{FP}^+(S_n)$ we can associate 
the binary relation
\begin{displaymath}
\varphi_A=\{(i,i\,\alpha):i\in N,\alpha\in A\}.
\end{displaymath}
It is easy to check that $\overline{A}\to\varphi_A$ is
a monomorphism from $\mathcal{FP}^+(S_n)$ to the semigroup
$B_n$ of all binary relations on $N$. Moreover,
for $A^{-1}:=\{\alpha^{-1}:\alpha\in A\}$ one has 
$\varphi_{A^{-1}}=\varphi_{A}^{-1}$. We will freely identify
$\mathcal{FP}^+(S_n)$ with its image in $B_n$ with respect
to this monomorphism.

The maximal nilpotent subsemigroups of $\mathcal{FP}^+(S_n)$
were described in \cite{GM2}. We will now recover this
description using Theorem~\ref{theorem4}.

In \cite[\S~3]{GM25} it is shown that the set of idempotents of 
$\mathcal{FP}^+(S_n)$ coincides with the set of all equivalence
relations on $N$. In particular, the zero element of 
$\mathcal{FP}^+(S_n)$ is the full relation $\theta=N\times N$.
Let $\rho$ and $\tau$ be two equivalence  relations.
It is easy to see that $\rho\tau=\tau\rho=\rho$ if
and only if each equivalence class of $\rho$ is the union of 
some equivalence classes of $\tau$. Hence, under the identification
of the equivalence relations with the correspondent decompositions
into equivalence classes, the natural partial order on idempotents
of the semigroup $\mathcal{FP}^+(S_n)$ coincides with the 
natural order in the lattice $\mathrm{Part}(N)$ of all
decompositions of the set $N$. Therefore the primitive idempotents
in $\mathcal{FP}^+(S_n)$ are those idempotents, which correspond to
decompositions of $N$ into $2$ blocks. It is obvious that all
primitive idempotents are minimal, and hence for 
the semigroup $\mathcal{FP}^+(S_n)$ these two notions coincide.
Let $\mathcal{M}$ denote the set of all
minimal (=primitive) idempotents of $\mathcal{FP}^+(S_n)$.

In \cite[Theorem~1]{Ma} it is shown that for arbitrary
$\alpha,\beta\in \mathcal{FP}^+(S_n)$ we have
$\alpha\mathcal{D}\beta$ is and only if 
$\alpha=\sigma\beta\delta$ for some $\sigma,\delta\in S_n$.
This immediately implies that two idempotents from
$\mathcal{FP}^+(S_n)$ are $\mathcal{D}$-related if and
only if the corresponding decompositions of $N$ have the
same type (that is the same number of subsets of each
cardinality). Hence the relation $\mathcal{D}$ induces 
the decomposition 
\begin{displaymath}
\mathcal{M}=\mathcal{M}_1\cup \mathcal{M}_2\cup\dots\cup
\mathcal{M}_{\lfloor\frac{n}{2}\rfloor},
\end{displaymath}
where $\mathcal{M}_k$, $1\leq k\leq \lfloor\frac{n}{2}\rfloor$,
consists of idempotents with two equivalence classes,
for which the ``smaller'' class contains exactly $k$ elements.
For $k<n/2$ we have $|\mathcal{M}_k|=\binom{n}{k}$, and
if $n=2l$ is even, we also have $|\mathcal{M}_l|=\frac{1}{2}\binom{n}{l}$.

For arbitrary $X\subset N$ and $\alpha\in \mathcal{FP}^+(S_n)$ set
\begin{displaymath}
X^{\alpha}=\{y:(x,y)\in\alpha\text{ for some }x\in X\}.
\end{displaymath}
Then for any $X\subset N$ and $\alpha\in \mathcal{FP}^+(S_n)$
we have $|X^{\alpha}|\geq |X|$ (see \cite[Lemma~8]{GM2}).
Moreover, if $|X^{\alpha}|= |X|$ for some $X\subset N$, then
$|(N\setminus X)^{\alpha}|=|N\setminus X|$ (see \cite[Lemma~1]{GM25}).

\begin{proposition}\label{prop12}
Let $\alpha\in \mathcal{FP}^+(S_n)$. Then
$\alpha\in R(\mathcal{FP}^+(S_n))$ if and only if
$|X^{\alpha}|>|X|$ for all proper subsets $X$ of $N$.
\end{proposition}

\begin{proof}
We start with necessity. Assume that there exists
a proper subset, $X$, of $N$ such that
$|X^{\alpha}|=|X|$. Then for the element
$\beta=\alpha\alpha^{-1}\in(\alpha)$ we have
\begin{equation}\label{eqfst}
X^{\beta}=X,\quad\quad 
(N\setminus X)^{\beta}=N\setminus X.
\end{equation}
This implies, in particular, that $X^{\beta^l}=X$ for every
$l\in\mathbb{N}$, which means that $\beta$ is not nilpotent,
and hence $(\alpha)$ is not a nilpotent ideal. This
proves the necessity.

To prove sufficiency we consider $\alpha\in \mathcal{FP}^+(S_n)$
such that $|X^{\alpha}|>|X|$ for all proper subsets $X$ of $N$.
Let $a\in N$. Consider arbitrary elements
$\mu_i\alpha\nu_i$, $i=1,\dots,n-1$, from $(\alpha)$. Then
for the element $\beta=\mu_1\alpha\nu_1\cdots \mu_{n-1}\alpha\nu_{n-1}$
we have 
\begin{multline}\label{eqfst.2}
1\leq |\{a\}^{\mu_1}|\lessdot|\{a\}^{\mu_1\alpha}|\leq
|\{a\}^{\mu_1\alpha\nu_1\mu_2}|\lessdot\\\lessdot
|\{a\}^{\mu_1\alpha\nu_1\mu_2\alpha}|\leq \dots\lessdot
|\{a\}^{\mu_1\alpha\nu_1\dots\mu_{n-1}\alpha}|\leq
|\{a\}^{\beta}|,
\end{multline}
where $x\lessdot y$ means $x=y$ if $x=n$ and $x<y$ if $x<n$.
In \eqref{eqfst.2} we have $n-1$ occurrences of $\lessdot$. Hence
$|\{a\}^{\beta}|\geq n$, which implies $\{a\}^{\beta}=N$ for
any $a\in N$. This means that $\beta=\theta$ and hence the
ideal $(\alpha)$ is a nilpotent ideal of nilpotency class
at most $n-1$. Therefore $\alpha\in R(\mathcal{FP}^+(S_n))$
and the proof is complete.
\end{proof}

Let us now show that the set $\mathcal{M}$ satisfies
the conditions \eqref{cond1}-\eqref{cond3}. 

To prove \eqref{cond1} let $\rho,\tau\in \mathcal{M}$,
$\rho\neq \tau$, be such that
$\rho$ corresponds to the decomposition $N=A_1\cup A_2$ and
$\tau$ to the decomposition $N=B_1\cup B_2$. If all four
intersections $A_i\cap B_i$ are non-empty, we obviously
have $\rho\tau=\tau\rho=\theta$. Assume now that
$A_1\cap B_1=\varnothing$. Then $A_1\subset B_2$,
$B_1\subset A_2$, $A_2=B_1\cup (A_2\cap B_2)$ and
$B_2=A_1\cup (A_2\cap B_2)$, moreover, $A_2\cap B_2\neq \varnothing$.
This implies that $\rho\tau=(N\times N)\setminus (A_1\times B_1)$.
Hence if $X\cap A_2=\varnothing$, we have $X^{\rho\tau}=N$.
If $X\subset A_1$, we have $X^{\rho\tau}=A_1\cup (A_2\cap B_2)$.
This means that $|X^{\rho\tau}|>|X|$ for any proper subset $X$ of $N$.
In particular, $\rho\tau\in R(\mathcal{FP}^+(S_n))$, and thus
$\rho$ and $\tau$ commute in the Rees quotient
$\mathcal{FP}^+(S_n)/R(\mathcal{FP}^+(S_n))$.

To prove \eqref{cond2} let $\alpha\not\in R(\mathcal{FP}^+(S_n))$.
According to Proposition~\ref{prop12}, there exist proper subsets
$X$ and $Y$ of $N$ such that $|X|=|Y|$, $X^{\alpha}=Y$ and
$(N\setminus X)^{\alpha}=N\setminus Y$. Consider the minimal idempotents
$\mu$ and $\eta$, which correspond to the decompositions
$N=X\cup (N\setminus X)$ and $N=Y\cup (N\setminus Y)$ respectively.
Obviously $X^{\mu\alpha\nu}=Y$ and hence $\mu\alpha\nu\not\in
R(\mathcal{FP}^+(S_n))$. Hence $\mu\alpha\nu\neq 0$ in the
Rees quotient $\mathcal{FP}^+(S_n)/R(\mathcal{FP}^+(S_n))$.

Finally, let us prove \eqref{cond3}. We first observe the
following: let $\mu$ be a minimal idempotent, which corresponds to the
decomposition $N=A_1\cup A_2$, and $X$ be a proper subset of
$N$ such that $|X^{\mu}|=|X|$;  then, obviously, $X=A_1$ or $X=A_2$.

Let now $\mu$ and $\nu$ be minimal idempotents, which correspond to
the decompositions $N=A_1\cup A_2$ and $N=B_1\cup B_2$
respectively. Let further $\alpha$ and $\beta$ be arbitrary elements
from $\mathcal{FP}^+(S_n)\setminus R(\mathcal{FP}^+(S_n))$.
Assume that $\mu\alpha\beta\nu\neq 0$ in
$\mathcal{FP}^+(S_n)/R(\mathcal{FP}^+(S_n))$. By Proposition~\ref{prop12},
there exists a proper subset, $X$, of $N$ such that
$|X^{\mu\alpha\beta\nu}|=|X|$. This means that we also have the 
following equalities:
\begin{displaymath}
|X|=|X^{\mu}|,\quad\text{ and }\quad
|(X^{\mu\alpha\beta})^{\nu}|=|X^{\mu\alpha\beta}|.
\end{displaymath}
In particular, we may assume $X=A_1$ and
$X^{\mu\alpha\beta}=B_1$. Therefore
$\mu\alpha\beta\nu\subset (A_1\times B_1)\cup(A_2\times B_2)$. Taking
into account that $\mu(\mu\alpha\beta\nu)=\mu\alpha\beta\nu$ and
$(\mu\alpha\beta\nu)\nu=\mu\alpha\beta\nu$ we even obtain
$\mu\alpha\beta\nu=(A_1\times B_1)\cup(A_2\times B_2)$.

Let now $\tau$ be the primitive idempotent, which corresponds
to the decomposition $A_1^{\alpha}\cup(N\setminus A_1^{\alpha})$.
Then by a direct calculation we have
\begin{displaymath}
\mu\alpha\tau\beta\nu=(A_1\times B_1)\cup(A_2\times B_2)=
\mu\alpha\beta\nu.
\end{displaymath}
Finally, let us assume that 
$\mu\alpha\rho\beta\nu\neq 0$ for some
primitive idempotent $\rho$. Then 
$|(A_i^{\alpha})^{\rho}|=|A_i^{\rho}|$ for $i=1,2$ and
hence $\rho=\tau$ and $\mu\alpha\rho\beta\nu=\mu\alpha\beta\nu$.
This proves that the condition \eqref{cond3} is satisfied.

Applying now Theorem~\ref{theorem4} and Corollary~\ref{cor6}
we obtain the following result:

\begin{theorem}\label{thm10}{\rm (\cite[Theorem~5]{GM2})}
There is a natural bijection between the maximal nilpotent
subsemigroups of the semigroup
$\mathcal{FP}^+(S_n)$ and the collections of linear
orders on the sets $\mathcal{M}_1,\dots,
\mathcal{M}_{\lfloor\frac{n}{2}\rfloor}$. In particular,
if $n=2k+1$, $\mathcal{FP}^+(S_n)$ contains
$\prod_{i=1}^k\binom{n}{i}!$ maximal nilpotent subsemigroups,
and if $n=2k$, $\mathcal{FP}^+(S_n)$ contains
$(\frac{1}{2}\binom{n}{k})!\prod_{i=1}^{k-1}\binom{n}{i}!$ 
maximal nilpotent subsemigroups.
\end{theorem}

\begin{remark}\label{remfps}
{\rm
Since $R(\mathcal{FP}^+(S_n))\neq \{0\}$, we are not able
to determine the nilpotency class of a maximal nilpotent subsemigroup
of $\mathcal{FP}^+(S_n)$ directly. However, this can be done by
different methods, see \cite{GM2} for details.
}
\end{remark}

\begin{remark}\label{remfps.2}
{\rm
The radical $R(\mathcal{FP}^+(S_n))$ appeared already in
\cite{GM2} as the intersection of all maximal nilpotent
subsemigroups of $\mathcal{FP}^+(S_n)$. It has some interesting
asymptotic properties, studied in \cite{GM3}.
}
\end{remark}

\subsection{Finite $0$-simple semigroups}\label{s6.10}

According to the celebrated theorem of Sushkevich and Rees
(see e.g. \cite[Theorem~5.3]{Gri}), a finite semigroup,
$S$, is $0$-simple if and only if it is isomorphic to some 
regular Rees matrix semigroup $M^{0}=M^{0}(G;I,\Lambda;P)$.
Moreover, the semigroup $M^{0}(G;I,\Lambda;P)$ is regular
if and only if the sandwich matrix $P$ is {\em regular}, that is
each row and each column of $P$ contains at least one
non-zero element.

All maximal nilpotent subsemigroups of the regular
semigroup $M^{0}$ were described in \cite{Gr} in two 
different ways: first in terms of some graphs, and
then in terms of some ordered partitions of the
index sets $I$ and $\Lambda$. The second description
reminds on the upper block-triangular matrices. In
\cite{Gr} one can also find an algorithm how to find
all necessary partitions of $I$ and $\Lambda$. Unfortunately,
neither the first not the second of these descriptions
allows one for example to estimate the number of 
maximal maximal nilpotent subsemigroups of $M^{0}$. 

It is obvious that the radical of a finite $0$-simple 
semigroup $S$ is $\{0\}$ and that every non-zero idempotent
of $S$ is minimal. However, in the general case we
can not apply Theorem~\ref{theorem5} since it is possible
that there is no subset of minimal idempotents,
satisfying \eqref{cond1}-\eqref{cond3}. Our first goal
is to give, in terms of the sandwich matrix $P$,
necessary and sufficient conditions for the
possibility of application of  Theorem~\ref{theorem5}
to the regular semigroup $M^{0}(G;I,\Lambda;P)$.

Let $I=\{1,2,\dots,n\}$ and $\Lambda=\{1',2',\dots,m'\}$,
$G$ be a finite group, and $P=(p_{j',i})$ be a regular sandwich 
matrix. Then the elements of $M^{0}(G;I,\Lambda;P)$ are triples
$(g,i,j')$, where $g\in G^0$, $i\in I$, and $j'\in\Lambda$
(we identify all triples of the form $(0,i,j')$). The
product of triples is defined as follows:
\begin{equation}\label{eqgra}
(g,i,j')(h,k,l')=(gp_{j',k}h,i,l').
\end{equation}

\begin{proposition}\label{prop7}
\begin{enumerate}[(a)]
\item\label{prop7.1}
The relation $\mathcal{H}$ on a regular Rees semigroup,
$M^{0}=M^{0}(G;I,\Lambda;P)$, is a congruence. The 
corresponding quotient semigroup $M^{0}/\mathcal{H}$ is
isomorphic to the regular Rees matrix semigroup
$\tilde{M}^{0}=$ $M^{0}(\{1\};$ $I,\Lambda;\tilde{P})$, where
$\{1\}$ is the group with one element, and
$\tilde{p}_{j',i}=1$ if and only if $p_{j',i}\neq 0$. 
\item\label{prop7.2}
The canonical epimorphism $\pi:M^{0}\to M^{0}/\mathcal{H}$
induces a bijection between the maximal nilpotent subsemigroups
of $M^{0}$ and the maximal nilpotent subsemigroups
of $\tilde{M}^{0}$. 
\item\label{prop7.3}
For every nilpotent subsemigroup $T\subset M^{0}$ the 
nilpotency classes of $T$ and $\pi(T)$ coincide.
\end{enumerate}
\end{proposition}

\begin{proof}
That $\mathcal{H}$ is a congruence on $M^{0}$ follows
from Lemma~\ref{lemma1} and Green's lemma (see e.g. 
\cite[Lemma~II.1.3]{Gri}). The isomorphism
$M^{0}/\mathcal{H}\cong \tilde{M}^{0}$ follows
from the definition \eqref{eqgra} of the multiplications
in $M^{0}$ and $\tilde{M}^{0}$, and the fact that 
the $\mathcal{H}$-classes of $M^0$ have the form
$H_{i,j'}=\{(g,i,j'):g\in G\}$ (see for example
\cite[Lemma~3.2]{CP}). This proves \eqref{prop7.1}.
\eqref{prop7.2} and \eqref{prop7.3} follow from the
fact that in $M^0$ we have $\mathcal{H}_0=\{0\}$.
\end{proof}

Proposition~\ref{prop7} implies, in particular, that every
maximal nilpotent subsemigroup of $M^0$ is a union of
$\mathcal{H}$-classes (this also follows from \cite[Theorem~1']{Gr}).
We will call the matrix $\tilde{P}$ from Proposition~\ref{prop7}
{\em reduced}. Since the first component of all non-zero elements
from $\tilde{M}^{0}$ equals $1$, we can consider these elements
as the pairs $(i,j')$, $i\in I$, $j'\in\Lambda$, with the following
modification of \eqref{eqgra}: $(i,j')(k,l')=(i,l')\tilde{p}_{j',k}$.

A $(0,1)$-matrix, $A$, will be called {\em row-regular} (resp.
{\em column-regular}) provided that each row (resp. column) of $A$
contains a non-zero element. Define the {\em Boolean product} of the
$(0,1)$-vectors $u=(u_1,\dots,u_k)$ and $v=(v_1,\dots,v_k)$
of the same length as follows:
\begin{displaymath}
u\circ v=
\begin{cases}
0, & u_1v_1+\dots+u_kv_k=0,\\
1, & \text{otherwise}.
\end{cases}
\end{displaymath}
Let now $A$ be a $(0,1)$-matrix of size $a\times b$ and
$B$ be a $(0,1)$-matrix of size $b\times c$. We define the 
{\em Boolean product} $A\circ B$ of $A$ and $B$ as the
$(0,1)$-matrix $C=(c_{i,j})$ of size $a\times c$, where
$c_{i,j}$ is the Boolean product of the $i$-th row of
$A$ with the $j$-th column of $B$. 

\begin{theorem}\label{theorem6}
The regular Rees matrix semigroup 
$\tilde{M}^{0}=(\{1\},I,\Lambda,\tilde{P})$ contains
a non-empty set of minimal idempotents, satisfying
the conditions \eqref{cond1}-\eqref{cond3}, if and only
if using independent permutations of rows and columns the
sandwich matrix $\tilde{P}$ can be reduced to the
form
\begin{equation}\label{eqmatrix}
\left(\begin{array}{c|c}
E_k& B \\\hline
A & A\circ B
\end{array}
\right),
\end{equation}
where $E_k$ is the identity matrix of size $k$, the matrix
$A$ is row-regular, and the matrix $B$ is column-regular.
If $\tilde{P}$ can be reduced to \eqref{eqmatrix}, then
$\tilde{M}^{0}$ contains exactly $k!$ maximal nilpotent
subsemigroups, each of which has nilpotency class $k$.
\end{theorem}

\begin{proof}
Assume that $\tilde{M}^{0}$ contains a set, $\{e_1,\dots,e_k\}$, 
of minimal idempotents, satisfying \eqref{cond1}-\eqref{cond3}.
From \cite[Lemma~3.2]{CP} we derive that the non-zero element
$(i,j')\in \tilde{M}^{0}$ is an idempotent if and only if
$\tilde{p}_{j',i}=1$. If $\tilde{p}_{j',i}=1$ then, abusing the
language, we will say that the idempotent $e=(i,j')$ is in
the $j'$-th row and $i$-th column of the matrix $\tilde{P}$.

According to Proposition~\ref{prop2}, minimal idempotents commute
if and only if they are orthogonal. Observe that the idempotents, 
which are in the same row (or in the same column) of $\tilde{P}$ 
can not commute. Indeed, if $e=(i,j')$ and $f=(i,k')$ are
idempotents, we have
\begin{displaymath}
ef=(i,j')(i,k')=(i,k')\tilde{p}_{j',i}=(i,k')\cdot 1=f\neq 0
\end{displaymath}
(analogously one shows that $e=(i,j')$ and $g=(k,j')$ do not
commute).

Permuting, if necessary, the elements of $I$ and $\Lambda$,
we can assume that $e_1=(1,1')$,\dots, $e_k=(k,k')$.

If $i\neq j$, then the commutativity of $e=(i,i')$ and
$f=(j,j')$ implies that $0=(i,i')(j,j')=(i,j')\tilde{p}_{i',j}$,
which yields $\tilde{p}_{i',j}=0$. Hence
$(\tilde{p}_{i',j})_{1\leq i,j\leq k}=E_k$ and
we have 
\begin{equation}\label{eqmatrix2}
\tilde{P}=\left(\begin{array}{c|c}
E_k& B \\\hline
A & C
\end{array}
\right)
\end{equation}
for some matrices $A$, $B$ and $C$. 

Let us now study the condition \eqref{cond2}.
For an arbitrary $(m,l')\in \tilde{M}^{0}$ there should exist
$e=(i,i')$ and $f=(j,j')$, $1\leq i,j\leq k$, such that
\begin{displaymath}
0\neq (i,i')(m,l')(j,j')=(i,j')\tilde{p}_{i',m}\tilde{p}_{l',j},
\end{displaymath}
that is there should exist $i'$ and $j$ such that 
$\tilde{p}_{i',k}=\tilde{p}_{l',j}=1$. The latter means that for 
arbitrary $m$ and $l'$ the $m$-th column of the matrix
$\left(E_k|B\right)$ and the $l'$-th row of the matrix
$\left(\begin{array}{c}E_k\\\hline A\end{array}\right)$ must 
be non-zero, that is the matrices $A$ and $B$ should be
row regular and column regular respectively.

Finally, let us look what does the condition \eqref{cond3} mean 
for the matrix \eqref{eqmatrix2}. Let $e=(i,i')$, $f=(j,j')$
and $g=(r,r')$, where $1\leq i,j,r\leq k$, $x=(u,v')$ and
$y=(p,q')$, where $u,v',p$ and $q$ are arbitrary. Then
\begin{displaymath}
exyf=(i,j')\tilde{p}_{i',u}\tilde{p}_{v',p}\tilde{p}_{q',j},
\quad\quad
exgyf=(i,j')\tilde{p}_{i',u}\tilde{p}_{v',r}\tilde{p}_{r',p}
\tilde{p}_{q',j}.
\end{displaymath}
If $\tilde{p}_{i',u}=0$ or $\tilde{p}_{q',j}=0$, then 
\eqref{cond3} is obviously satisfied. Let now
$\tilde{p}_{i',u}=\tilde{p}_{q',j}=1$. If $\tilde{p}_{v',p}=0$
then for each $r$, $1\leq r\leq k$, we should have
$\tilde{p}_{v',r}\tilde{p}_{r',p}=0$. Hence the Boolean product
of the $v'$-th row of the matrix
$\left(\begin{array}{c}E_K\\\hline A\end{array}\right)$
with the $p$-th column of the matrix $\left(E_k|B\right)$ should
be $0$, that is $\tilde{p}_{v',p}$. This is obvious if
$1'\leq v'\leq k'$ or $1\leq p\leq k$. If both
$v'>k'$ and $p>k$, then $\tilde{p}_{v',p}$ is a coefficient of
the matrix $C$ and we get the restriction that  this coefficient must
be the Boolean product of the corresponding row of $A$ with the
corresponding column of $B$.

If $\tilde{p}_{v',p}=1$, there should exist some $r$, 
$1\leq r\leq k$, such that $\tilde{p}_{v',r}\tilde{p}_{r',p}=1$.
This is again obvious if $1'\leq v'\leq k'$ or $1\leq p\leq k$. If both
$v'>k'$ and $p>k$, then $\tilde{p}_{v',p}$ is again a coefficient of
the matrix $C$ and we again get the restriction that this coefficient must
be the Boolean product of the corresponding row of $A$ with the
corresponding column of $B$. Hence $C=A\circ B$ and the necessity is 
proved. 

To prove sufficiency one checks by a direct calculation that,
given \eqref{eqmatrix}, the idempotents $(1,1')$,\dots, $(k,k')$
satisfy \eqref{cond1}-\eqref{cond3}. 

The rest follows now immediately from Theorem~\ref{theorem5}.
\end{proof}

If the sandwich matrix $\tilde{P}$ has the form \eqref{eqmatrix} then
one can explicitly describe all elements of the maximal nilpotent 
subsemigroup, which correspond to the natural linear order
$(1,1')<(2,2')<\dots<(k,k')$ on our set of minimal idempotents 
(note that all other linear orders can be
reduced to this one by permutations of rows and columns). For this
we consider the new matrix
\begin{equation}\label{eqmatrix3}
\tilde{P}^*=(\tilde{p}^*_{j',i})=\left(\begin{array}{c|c}
UT_k& B^* \\\hline
A^* & A^*\circ B^*
\end{array}
\right),
\end{equation}
which is obtained from \eqref{eqmatrix} in the following way:
in every row of $A$ we look for the leftmost occurrence of $1$ and all
the $0$'s to the right of this $1$ we change to $1$'s; in 
every column of $B$ we look for the lowest occurrence of $1$ and
all the $0$'s above this $1$ we change to $1$'s; the matrix $UT_k$
is the matrix of order $k$ which has $0$'s below the main diagonal
and $1$'s everywhere else. Observe that $UT_k$ is obtained from
$E_k$ by the same rule which we have just used to create $A^*$ 
from $A$ (or $B^*$ from $B$).

\begin{proposition}\label{prop8}
Let $T$ be the maximal nilpotent subsemigroup of
$\tilde{M}^{0}$, which corresponds 
to the natural linear order $(1,1')<(2,2')<\dots<(k,k')$
on the set $\{(1,1'),\dots,(k,k')\}$ of minimal idempotents.
Then a non-zero element, $x=(u,v')\in\tilde{M}^0$, belongs to 
$T$ if and only if $\tilde{p}^*_{v',u}=0$.
\end{proposition}

\begin{proof}
First we observe that
\begin{displaymath}
\tilde{P}^*=
\left(\begin{array}{c}UT_k\\\hline A^*\end{array}\right)\circ
\left(UT_k|B^*\right).
\end{displaymath}
By Theorem~\ref{theorem5}, $x\in T$ if and only if we have
\begin{displaymath}
(i,i')(u,v')(j,j')=(i,j')\tilde{p}_{i',u}\tilde{p}_{v',j}\neq 0
\end{displaymath}
implies $i<j$ for all $1\leq i,j\leq k$. Let now $i_1'$ be the
maximal element $i'$ of $\{1',2',\dots,k'\}$ such that 
$\tilde{p}_{i',u}=1$, and $j_1$ be the minimal element $j$ of
$\{1,2,\dots,k\}$ such that $\tilde{p}_{v',j}=1$.
As $\tilde{p}_{i_1',u}\tilde{p}_{v',j_1}=1$, we have that
$x\in T$ if and only if $i_1<j_1$, that is if and only if the Boolean
product $\tilde{p}^*_{v',u}$ of the $v'$-th row of the matrix
$\left(\begin{array}{c}UT_k\\\hline A^*\end{array}\right)$ with
the $u$-th column of the matrix $\left(UT_k|B^*\right)$ equals
$0$.
\end{proof}

\begin{corollary}\label{cor10}
Assume that the reduced sandwich-matrix of the regular semigroup
$M^0=M^0(G,I,\Lambda,P)$ has the form \eqref{eqmatrix}. Let
$e_i$, $1\leq i\leq k$, denote the idempotent from the class
$\mathcal{H}_{i,i'}$. Let further $T$ be the maximal nilpotent
subsemigroup of $M^0$, which corresponds to the linear order
$e_1<e_2<\dots<e_k$. Then
\begin{displaymath}
T=\{0\}\cup\bigcup_{\tilde{p}^*_{v',u}=0}\mathcal{H}_{u,v'}.
\end{displaymath}
In particular, $|T|=1+f\cdot|G|$, where $f$ denotes the number
of $0$'s in the matrix $\tilde{P}^*$.
\end{corollary}

\begin{remark}\label{remlast}
{\rm
Note that for $k=1$ the matrix \eqref{eqmatrix} does not
contain any zero entries, and for $k>1$ it contains at least $2$
zero entries. Hence in the case when the sandwich-matrix of the 
regular semigroup $M^0=M^0(G,I,\Lambda,P)$ contains exactly one
zero entry, then the reduced sandwich matrix can not be
written in the form \eqref{eqmatrix}. The same is also the case if
the original sandwich matrix contains more than two $0$'s which
are all in different rows (or columns), or in the same row
(or column).
}
\end{remark}

\subsection{The semigroup of square matrices}\label{s6.11}

Consider now the semigroup $\mathrm{Mat}_n(\mathbb{F}_q)$ of
all $n\times n$ matrices over the field $\mathbb{F}_q$ with
$q<\infty$ elements. Maximal nilpotent subsemigroups of
$\mathrm{Mat}_n(\mathbb{F}_q)$ were classified in
\cite[Corollary~3]{KM}. It turns out that there is a bijection
between maximal nilpotent subsemigroups of
$\mathrm{Mat}_n(\mathbb{F}_q)$ and complete flags in
the vector space $\mathbb{F}_q^n$. In particular,
$\mathrm{Mat}_2(\mathbb{F}_2)$ contains $3$ maximal
nilpotent subsemigroups. Since, obviously,
$R(\mathrm{Mat}_2(\mathbb{F}_2))=\{0\}$ and $3$ is not a 
product of factorials, from  Theorem~\ref{theorem5} we deduce 
that there does not exist any subset of the set of minimal 
idempotents in $\mathrm{Mat}_2(\mathbb{F}_2)$, which satisfies 
\eqref{cond1}-\eqref{cond3} (and then this kind of argument can 
be extended to an arbitrary $\mathrm{Mat}_n(\mathbb{F}_q)$). 
On the other hand, it
is easy to find a set of minimal idempotents in
$\mathrm{Mat}_n(\mathbb{F}_q)$, which satisfies
both \eqref{cond1} and \eqref{cond2}. For $1\leq i,j\leq n$
let $e_{i,j}$ denote the corresponding {\em matrix unit} 
(that is the matrix which has only one non-zero
entry, namely $1$, staying on the intersection of $i$-th row and
$j$-th column). Then it is easy to see that the set 
$\{e_{1,1},\dots,e_{n,n}\}$ satisfies both \eqref{cond1} and 
\eqref{cond2}. In particular, this shows that the condition
\eqref{cond3} is really important for Theorem~\ref{theorem5}.


\vspace{1cm}

\noindent
O.G.: Department of Mechanics and Mathematics, Kyiv Taras Shevchenko
University, 64, Volodymyrska st., 01033, Kyiv, UKRAINE,\\
e-mail: {\tt ganiyshk\symbol{64}univ.kiev.ua}
\vspace{0.5cm}

\noindent
V.M.: Department of Mathematics, Uppsala University, Box 480,
SE 751 06, Uppsala, SWEDEN, e-mail: {\tt mazor\symbol{64}math.uu.se},\\
web: ``http://www.math.uu.se/$\tilde{\hspace{2mm}}$mazor''
\vspace{0.5cm}

\end{document}